\documentclass[11pt,twoside,onecolumn]{article}

\usepackage{authblk} 

\usepackage{graphicx}
\graphicspath{ {Vis/} }
\usepackage[hang, small,labelfont=bf,up,textfont=it,up]{caption}
\usepackage{subcaption,dsfont}

\usepackage[usenames, dvipsnames]{color}

\usepackage{indentfirst} 

\usepackage[sc]{mathpazo} 
\usepackage[T1]{fontenc} 
\usepackage{microtype} 

\usepackage[english]{babel} 

\usepackage[hmarginratio=1:1,top=32mm,columnsep=20pt]{geometry} 
\usepackage{booktabs} 

\usepackage{lettrine} 

\usepackage{enumitem} 
\setlist[itemize]{noitemsep} 

\usepackage{titlesec} 

\usepackage{fancyhdr} 
\usepackage{hyperref} 

\usepackage{mathtools} 

\usepackage{amsmath,calc,amsthm}

\DeclareMathOperator*{\argmin}{argmin}

\usepackage{amsfonts}
\usepackage{amsthm}

\theoremstyle{definition}

\usepackage{multicol}
\usepackage{multirow}
\usepackage{amssymb}
\usepackage{array}
\usepackage{tabu}
\usepackage{algpseudocode}
\usepackage[numbers]{natbib}

\theoremstyle{definition}

\theoremstyle{plain}

\usepackage{xcolor}

\usepackage{makecell}
\usepackage{scrextend} 

\usepackage{algorithm}
\usepackage{nomencl}
\usepackage{etoolbox}
\makenomenclature

\makeatletter
\patchcmd{\thenomenclature}{\section*}{\section}{}{}
\makeatother

\begin{document}

\title{Approximate Dynamic Programming\\
	   for Real-time Dispatching and Relocation\\
	   of Emergency Service Engineers} 

\author[1]{D. Usanov\thanks{Corresponding author. Tel.:+31(0)20 592 4168.\\
\textit{E-mail addresses}: usanov@cwi.nl (D. Usanov), pechina.anna@gmail.com (A. Pechina), p.m.van.de.ven@cwi.nl (P.M. van de Ven), r.d.van.der.mei@cwi.nl (R.D. van der Mei).}}
\author[1]{A. Pechina}
\author[1]{P.M. van de Ven}
\author[1, 2]{R.D. van der Mei}
\affil[1]{Centrum Wiskunde \& Informatica, Science Park 123,
1098 XG Amsterdam, The Netherlands}
\affil[2]{Vrije Universiteit Amsterdam, De Boelelaan 1081a,
1081 HV Amsterdam, The Netherlands }

\date{} 


\maketitle
\newpage
\begin{abstract}
Quick response times are paramount for minimizing downtime in spare parts networks for capital goods, such as medical and manufacturing equipment. To guarantee that the maintenance is performed in a timely fashion, strategic management of both spare parts and service engineers is essential. While there is a rich body of research literature devoted to spare parts management, the problem of real-time management of service engineers has drawn relatively little attention. Motivated by this, we consider how to dispatch service engineers to breakdowns, and how to relocate idle engineers between base stations. We develop an approximate dynamic programming (ADP) approach to produce dispatching and relocation policies, and propose two new algorithms to tune the ADP policy. We conduct extensive computational experiments to compare the ADP policy against two benchmark policies by means of simulation. These demonstrate that the ADP approach can generate high-quality solutions that outperform both benchmarks across a wide range of networks and parameters. We observe significant improvements in terms of fraction of late arrivals over the two benchmarks, without increase in average response time.

\end{abstract}



\section{Introduction}
Manufacturers of capital goods typically provide post-sale support, in order to avoid frequent downtime of these expensive and essential goods. For instance, manufacturing plants that are down may result in significant costs for the owner, and an MRI scanner that is not working may even endanger lives. Such post-sale support accounts for  significant fraction of the profits of capital-goods manufacturers~\cite{koudal2006}.

In order to ensure minimal downtime, manufacturers often provide corrective maintenance, repairing and replacing machine parts that are malfunctioning. To meet the service-level agreements with their clients, these manufacturers must be able to quickly dispatch a service engineer and the necessary spare parts upon being informed of a breakdown. Doing so requires an extensive network of service engineers and spare part supply depots, as well as clever policies for managing these. This entails striking a careful balance between meeting the agreements and maintaining low costs~\cite{murthy2004}.

Designing and managing the spare parts side of this problem has been extensively studied in the research literature (cf.~\cite{houtum2015}), but the question of how to deal with service engineers in a cost-effective manner has received relatively little attention. These engineers are positioned in geographically dispersed locations, such that they jointly cover all possible locations of a breakdown. However, when one or more service engineers are busy doing maintenance, gaps in the coverage may arise, leading to costs incurred for missing service level agreements. We are interested in the question of how to manage service engineers in the most efficient way. This includes decisions on which engineer to dispatch to a new maintenance job, whether to relocate engineers to different base stations when gaps in the coverage arise, and whether to hold off on certain demand in order to wait for nearby service engineers to complete their existing job.

To our knowledge, the first such model was introduced in~\cite{drent2018}, where the authors studied a grid-like type of network and developed an algorithm for dynamic movement of engineers on the grid that outperformed the static closest-first policy. These kinds of decision problems can be written as a Markov decision process (MDP) in order to find the optimal dispatching and relocation policy, but in practice the so-called curse of dimensionality makes this approach often infeasible for real life-sized problems. In~\cite{pechina2019} the authors consider a general network structure and compare a wide range of heuristics taken from the research on a closely related problem of dispatching and relocation of Emergency Medical Services (EMS)~\cite{pechina2019}, one showed the best performance in almost all cases. In our work we use that algorithm as one of the benchmark policies.

The contribution of this paper is threefold. First, we develop a new relocation and dispatching policy for the model studied in~\cite{pechina2019}. The policy is based on Approximate Dynamic Programming (ADP), which relies on approximating the value function found in MDPs by a linear combination of easily computable basis functions. When the basis functions and their coefficients are selected carefully, such an approximation may yield sub-optimal but excellent policies that do not suffer from the scalability issues encountered with MDPs. For an excellent discussion of various ADP techniques, we refer to~\cite{powell2007approximate}. Second, we introduce two algorithms for fine-tuning ADP in our setting. The two algorithms, genetic algorithm and tabu search, tune the coefficients of the basis functions based on the actual performance of the corresponding policy, rather than on making a close approximation of the value function. Finally, we conduct extensive computational experiments, where the ADP policies obtained with both tuning algorithms are compared against two benchmark policies for various types of systems of realistic size. The benchmarks policies are the heuristic algorithm from~\cite{pechina2019} and the closest-first dispatching policy commonly used both in practice and, as a benchmark policy, in research literature. We show that the ADP policies outperform both benchmarks on various types of networks. By tuning the ADP cost function appropriately it is possible to significantly reduce the fraction of late arrivals while maintaining the same level of average response time observed under the benchmark policies. 

The remainder of the paper is organized as follows. Section \ref{sec:lit} gives an overview of related literature. In Section \ref{sec:model}, we present the model and formulate it as a Markov decision process. The benchmark heuristic is described in Section~\ref{sec:heuristic}. Section \ref{sec:adp} introduces the ADP approach together with the tuning algorithms. Numerical experiments are presented in Section \ref{sec:num}. Finally, Section \ref{sec:conclusion} contains concluding remarks and discussion.

\section{Literature review}\label{sec:lit}

There are three streams of literature most related to the work in this paper. The one that is contextually closest is the stream of research in the area of spare parts management. For an overview of this area we refer to~\cite{houtum2015}. Some studies in this stream that are most relevant to our work are those optimizing lateral transshipment policies. Lateral transshipment is a stock transfer from one local warehouse to another one. Examples of such studies, focusing on how to make transshipment decisions, include~\cite{axsater2003, minner2003, lee2007}. Another closely related model was recently proposed in~\cite{tiemessen2013}. The authors studied inventory networks with multiple customer classes where demand can be fulfilled directly (rather than via lateral transshipment).

The second stream of literature is in the field of dynamic EMS management. It is well-studied and is closely related to our setting. Therefore, we outline recent developments in the area of EMS management. For an extensive overview of literature on location, relocation and dispatching of EMS, we refer to~\cite{ingolfsson2013, belanger2019}. A common way to model EMS systems used in literature is MDP~\cite{mclay2013b, bandara2014, jagtenberg2017a}. In~\cite{bandara2014} and \cite{jagtenberg2017a}, for example, the authors demonstrate that the commonly used in practice policy of always dispatching the closest ambulance is not necessarily optimal. Although insightful, using MDP to find optimal policies is typically computationally intractable for real world applications. Therefore, researchers often resort to heuristic algorithms to make relocation and dispatching decisions. In~\cite{gendreau2006}, the so-called Maximum Expected Coverage Relocation Problem (MEXCRP) was formulated to compute compliance tables for ambulance relocation. A compliance table precomputes decirable vehicle locations depending on the number of idle emergency vehicles in the system. The relocations are then made according to the compliance table whenever this number changes, e.g., when a vehicle is dispatched to an incident or finishes a job. The algorithm was later extended in~\cite{barneveld2016b} to the Adjusted Maximal Expected Coverage Relocation Problem (AMEXPREP) by incorporating the busy fraction of the ambulances.

Compliance tables are easy to use in practice, as those are computed only once, and no extra computations are typically needed when applying the obtained policy in real time. The disadvantage of such policies however is that they incorporate a limited amount of information about the current state of the system, ignoring, for example, the spacial distribution of idle vehicles. Alternative approaches aim at deriving decisions in real time for a given state of the system. One such real-time relocation model was introduced in~\cite{gendreau2001}, that used the Double Standard Model~\cite{gendreau1997} to maximize the demand covered by at least two vehicles when making relocations while minimizing the relocation costs. In~\cite{andersson2007}, the authors introduced the notion of preparedness, that is, a measure characterizing the ability of the system to respond to current and future incidents. They introduced the dispatching and relocation algorithms maximizing the preparedness of the system.

In~\cite{jagtenberg2015}, the Dynamic Maximal Expected Coverage Location Problem (DMEXCLP) was introduced for redeployment of ambulances that just finished their job. The heuristic makes decisions that lead to a better configuration of the system in terms of expected covered demand. Later, in~\cite{barneveld2016a} the algorithm was further extended to incorporate the opportunity for a relocation right after a dispatching decision is made. Similar logic was used in~\cite{jagtenberg2017b} for making dispatching decisions, where each candidate ambulance for dispatching was evaluated in terms of the remaining expected coverage without it. An extension of this model in~\cite{barneveld2017} also included the possibility for postponing a dispatching decision until a closer busy ambulance finishes its current job.

Over the last decade, the ADP approach was succesfully applied to the problem of ambulance dispatching and relocation. In~\cite{maxwell2010}, ADP was used for ambulance redeployment upon completion of service, assuming a fixed closest-first dispatching policy and no other relocations. In~\cite{schmid2012}, the authors used ADP to optimize both dispatching and redeployment of ambulances that finish their jobs. Lately, in~\cite{nasrollahzadeh2018} the authors formulated a general model where the new incidents can be put in a queue, and there is a possibility to relocate idle ambulances.

The literature on ADP is the third stream of research most related to our work. Apart from the above mentioned applications in dynamic EMS management, ADP was successfully used in other application areas. Some important examples include fleet management~\citep{simao2009b}, dynamic container allocation~\citep{lam2007}, spare parts and supply chain management~\citep{simao2009a, fang2013}, capacity allocation in service industry~\citep{schuetz2012} and healthcare~\citep{astaraky2015}.

\section{Model description}\label{sec:model}
In this section we recapitulate the model introduced in \citep{pechina2019}. We consider a service region consisting of a set of identical machines $\mathcal{K} = \{1,...,K\}$ and a set of base stations $\mathcal{R} = \{1,...,R\}$. Locations of the machines and the base stations are fixed, and the traveling times between each pair of locations are deterministic and known. We assume that each machine has at least one base station that is reachable within the given time limit $TL$ (i.e., each machine is covered by at least one station). We also assume that each base station covers at least one machine. Let $\mathcal{M} = \{1,...,M\}$ denote the set of service engineers. Each service engineer is rested at one of the base station when idle.

The time till the next breakdown of a working machine is exponentially distributed with rate $\lambda$. Upon a breakdown of a machine, exactly one service engineer is needed to repair it, and the repair time of each machine is exponentially distributed with rate $\mu$. This repair time does not include the traveling time required for a service engineer to reach the location of a machine. However, the repair starts immediately after a service engineer arrives to the location of a broken machine, assuming all the necessary tools and spare parts needed for repair are available upon arrival. A service engineer can be dispatched immediately after a failure occurred, or the failed machine can be put in a queue waiting for service. The time between the failure and the arrival of a service engineer to the failure location is called \textit{response time}. If the response time exceeds the time limit $TL$, costs are incurred. We discuss the costs structure later in this section.

We consider the state of the system immediately after one of the following events happen:
\begin{enumerate}
\item a failure of a machine;
\item end of repair;
\item arrival of a service engineer at a base station;
\item arrival of a service engineer at a machine location.
\end{enumerate}
The event $e$ is described by the tuple $(e_t, e_l)$, where $e_t \in \{1, 2, 3, 4\}$ indicates the type of the event, and $e_l \in \{1,...,L\}$ the event location. The location can be either a broken machine or a base station, with the total number of locations $L = K+R$. Let $l_k$ indicate the location of machine $k \in \mathcal{K}$, and $l_r$ the location of base station $r \in \mathcal{R}$.

Let $\kappa_k$ denote the state of machine $k\in \mathcal{K}$. Here $\kappa_k = 0$ if machine $k$ is working, $\kappa_k = -1$ if machine $k$ is in repair, and $\kappa_k = t$ if machine $k$ has been waiting for repair for $t$ time units. Note that immediately after failure of machine $k$ its state is $\kappa_k = 0$ despite the machine is broken, as its elapsed waiting time is $0$ time units. This plays a role when defining the number of broken/working machines in a given state (i.e., if the event type is $e_t = 1$, the number of broken machines is $\big|\big\{ k \in \mathcal{K}: (\kappa_k  \neq 0) \lor (e_l = l_k)  \big\}\big|$). The state of all machines is represented by the vector $\pmb{\kappa} = (\kappa_1,...,\kappa_K)$.

We describe the state of all service engineers by the vector $\pmb{\mathfrak{m}} = (\mathfrak{m}_1,..,\mathfrak{m}_M)$. Here $\mathfrak{m}_m = (l_m, d_m)$ describes the state of service engineer $m$, where $l_m$ indicates the destination of that service engineer, and $d_m$ indicates the remaining time left to reach the destination. If the service engineer $m$ is residing at a location $l$ (either doing a repair or waiting at a base station), then $l_m = l$ and $d_m = 0$.

The state of the system immediately after an event $e$ is described by the tuple $(t, e, \pmb{\kappa}, \pmb{\mathfrak{m}})$, where $t$ is the time of the event. We use $t(s)$, $e(s)$, $\pmb{\kappa}(s)$, $\pmb{\mathfrak{m}}(s)$ to represent each component of a given state $s$. The state space is infinitely large, with non-recurrent states, as the time is included in the state description.

\subsection{Actions}\label{sec:actions}
Whenever an event occurs, an action is taken. The set of possible actions is defined by the type of event $e_t$. Below we formulate the action space per event type, but first we introduce additional notation and general assumptions regarding feasible actions. We assume that only idle service engineers are allowed to be dispatched to a broken machine. We denote the set of all idle service engineers in state $s$ by $\mathcal{F}(s) = \big\{m \in \mathcal{M} \mid l_m(s) \in \mathcal{R}\big\}$. Note that the service engineer who is still on the way to the base station is considered idle, and can be dispatched for repair. In that case, we assume that the service engineer first reaches the base station and then immediately departs to the machine location. We assume, however, that relocation of traveling idle service engineers is not allowed. We do this to avoid situations when an idle service engineers is continuously relocated from one station to another, never reaching his/her destination. Let $\mathcal{F}_0(s) = \big\{m \in \mathcal{M} \mid l_m(s) \in \mathcal{R}, d_m = 0\big\}$ denote the subset of idle service engineers that are not traveling. Let also denote the set of broken machines that are waiting in the queue in state $s$ by $\mathcal{Q}(s) = \big\{k \in \mathcal{K} \mid \{m \in \mathcal{M} \mid l_m(s)=k\} = \emptyset\big\}$.\\

\noindent
\textbf{Event type \pmb{$e_t = 1$}}. In case of a new failure, the action consists of a dispatching decision and a relocation decision. At this point, dispatching is allowed only to the newly broken machine, but not to the machines in the queue. One of the idle service engineers $\mathcal{F}(s)$ can be dispatched to the new breakdown, or a machine can be put in the queue. If dispatching is done, a relocation is allowed of one of the remaining stationary idle service engineers from his/her current location to another base station. However, if the repair request is put in the queue, relocation is not allowed. Let the binary vector $\pmb{X}$ of length $M$ represent the dispatching decision, and the binary $M \times R$ matrix $\pmb{Y}$ the relocation decision. Here $X_m = 1$ indicates that service engineer $m$ is dispatched for repair, and $Y_{mr} = 1$ indicates that service engineer $m$ is relocated to base station $r$. The action space can be formally described as
\begin{align}
\nonumber \mathcal{A}_1(s) = \bigg\{(\pmb{X}, \pmb{Y}) \mid &X_m = 0, \ \ m \notin \mathcal{F}(s); \ Y_{mr} = 0, \ \ m \notin \mathcal{F}_0(s),\ r \in \mathcal{R};  \\
\nonumber &\sum_{m \in \mathcal{M}}X_m \le 1; \ \sum_{\substack{m \in \mathcal{M},\\ r \in \mathcal{R}}} Y_{mr}\le 1;  \\ 
&\sum_{\substack{m \in \mathcal{M},\\ r \in \mathcal{R}}} X_mY_{mr}=0; \ \Big(1-\sum_{m \in \mathcal{M}}X_m\Big)\sum_{\substack{m \in \mathcal{M},\\ r \in \mathcal{R}}} Y_{mr} = 0 \bigg\}, \label{eqn:A1}
\end{align}
where the constraints ensure that at most one idle service engineer is dispatched, at most one other idle stationary service engineer is relocated, and relocation is done only upon dispatching.\\

\noindent
\textbf{Event type \pmb{$e_t = 2$}}. When a service engineer finishes repairing a machine, we can either allocate that service engineer to one of the base stations or dispatch to one of the broken machines in the queue. No relocation of other service engineers is allowed in this case. In \citep{pechina2019}, the authors included the possibility of relocating one other idle service engineer for this type of event. However, the best performing heuristic in that study did not make use of such extra relocations. Assuming that the potential gain from one extra relocation is marginal, we exclude the possibility of extra relocation to reduce computational complexity of the ADP approach. Moreover, using the same action space for both approaches makes the comparison cleaner.

Let the binary vector $\pmb{Z}$ of length $L$ represent the redeployment decision. Here $Z_l = 1$ indicates that the service engineer that just finished a repair is redeployed to location $l$, where $l$ is the location of either one of the base stations or one of the machines in the queue. The corresponding action space is given by
\begin{equation}\label{eqn:A2}
\mathcal{A}_2(s) = \bigg\{\pmb{Z} \mid  Z_l = 0, \ \ l \notin \{l_k: k \in \mathcal{Q}(s)\} \cup \{l_r: r \in \mathcal{R}\}; \sum_{l \in \mathcal{L}}Z_l = 1 \bigg\},
\end{equation}
where the constraints ensure that the service engineer is redeployed to either a base station or a machine from the queue.\\

\noindent
\textbf{Event type \pmb{$e_t = 3$}}. When a service engineer arrives at a base station, he/she can be immediately dispatched to one of the machines in the queue or left idle at that base station. Relocation of any sort is not allowed in that case. Let the binary vector $\pmb{U}$ of length $K$ represent the dispatching decision. Here $U_k=1$ indicates that the service engineer that just arrived at a base station is dispatched to machine $k \in \mathcal{Q}(s)$. The action space can be written as
\begin{equation}\label{eqn:A3}
\mathcal{A}_3(s) = \bigg\{ \pmb{U} \mid U_k = 0, \ \ k \notin \mathcal{Q}(s);\sum_{k \in \mathcal{K}}U_k \le 1 \bigg\},
\end{equation}
where the constraints ensure that the service engineer is dispatched to at most one machine from the queue. Note that $\mathcal{A}_3(s)$ is empty if the queue is empty in state $s$.\\

\noindent
\textbf{Event type \pmb{$e_t = 4$}}. Once a service engineer arrives at a broken machine, he/she immediately starts repairing the machine. In that case no action is taken, hence
\begin{equation}\label{eqn:A3}
\mathcal{A}_4(s) = \emptyset.
\end{equation}

\subsection{Transitions}\label{sec:trans}
We describe evolution of the system via the discrete-time stochastic process $\{s_n\}_{n\in \mathds{N}}$ embedded on decision epochs. To that end, we need to derive the distribution of time till the next transition for a given state $s_n$ and action $a_n$, as well as the transition probabilities.

Given state $s_n$ and action $a_n$, the next state $s_{n+1}$ is defined by the function $s_{n+1} = \Phi(s_n, a_n, \omega(s_n, a_n))$, where the random element $\omega(s_n, a_n)$ determines the next event. The next event can be either a failure of one of the working machines, the end of an ongoing repair, or an arrival of one of the traveling service engineers to his/her destination. Let $d(s_n, a_n)$ denote the minimum remaining distance in time units over all traveling service engineers after action $a_n$ is taken in state $s_n$. Let $d(s_n, a_n) = \infty$ if there are no traveling service engineers. Let also $\mathcal{W}(s_n)$ denote the set of all working machines and $\mathcal{H}(s_n)$ the set of machines in repair in state $s_n$. Note that those sets are not affected by the action $a_n$. Recall that the time till breakdown of each machine in $\mathcal{W}(s_n)$ is exponentially distributed with rate $\lambda$, and the time till the end of repair of each machine in $\mathcal{H}(s_n)$ is exponentially distributed with rate $\mu$. Hence, the time till the next event is distributed as the minimum of $d(s_n, a_n)$ and an exponentially distributed random variable $\Gamma(s_n)$ with rate $\eta(s_n) = \lambda|\mathcal{W}(s_n)| + \mu|\mathcal{H}(s_n)|$.

The probability that the next event is of type $e_t(s_{n+1}) = 1$ or $e_t(s_{n+1}) = 2$ is equal to $P(\Gamma(s_n)<d(s_n, a_n)) = 1-e^{-\eta(s_n) d(s_n,a_n)}$. Then the probability that the next event is
\begin{itemize}
\item the failure of machine $k \in \mathcal{W}(s_n)$ is $$\frac{\lambda}{\eta(s_n)}(1-e^{-\eta(s_n) d(s_n,a_n)}),$$ if $\mathcal{W}(s_n)\ne\emptyset$, and $0$ otherwise;
\item the end of repair of machine $k \in \mathcal{H}(s_n)$ is $$\frac{\mu}{\eta(s_n)}(1-e^{-\eta(s_n) d(s_n,a_n)}),$$ if $\mathcal{H}(s_n)\ne\emptyset$, and $0$ otherwise;
\item the arrival of a service engineer to his/her destination is $$e^{-\eta(s_n) d(s_n,a_n)},$$ if there are traveling service engineers after action $a_n$ is taken in state $s_n$, and $0$ otherwise.
\end{itemize}

The state of all service engineers immediately after action $a_n$ is taken in state $s_n$ is updated in a straightforward manner by changing destinations and distances accordingly. Upon realization of the next event $e(s_{n+1})$ and transition time $\min\{\Gamma(s_n), d(s_n,a_n)\}$, the state $s_{n+1}$ is obtained as follows. The time is equal to $t(s_{n+1}) = t(s_n) + \min\{\gamma(s_n), d(s_n,a_n)\}$, where $\gamma(s_n)$ is the realization of $\Gamma(s_n)$. If the next event is the repair of machine $k$, its state is set to $\kappa_k = 0$. Waiting times of all machines in the queue are increased by $\min\{\gamma(s_n), d(s_n,a_n)\}$, and traveling times of all traveling service engineers are decreased by $\min\{\gamma(s_n), d(s_n,a_n)\}$. If the next even is an arrival of a service engineer at a machine $k$, the state of that machine is set to $\kappa_k = -1$.

\subsection{Costs} \label{sec:costs}
Our main objective is to maximize the fraction of failures responded to within the time limit $TL$. We define the cost structure accordingly. If a machine breaks down, and a service engineer does not reach this machine within the time limit $TL$, then a penalty $1$ is paid. Apart from that, a small penalty $0<\epsilon \ll 1$ is paid per time unit of service delay over $TL$. The penalty $\epsilon$ is introduced only to prevent the situation of leaving some machines broken for a very long time or even forever, which is optimal in long-term perspective in terms of fraction of late arrivals but not realistic. So $\epsilon$ should not be too small to avoid unreasonably large waiting times.

Let $c\left(s_n, a_n, s_{n+1}\right)$ denote the costs incurred during the transition period from state $s_n$ to state $s_{n+1}$ when action $a_n$ is taken. The transition costs are computed as follows. A unit cost is incurred for each broken machine whose waiting time exceeds the time limit $TL$ during the transition. In addition, extra penalty $\epsilon$ is incurred per unit of waiting time over $TL$ for each broken machine. Then in total
\begin{equation}\label{eq:costs}
c(s_n, a_n, s_{n+1}) = \sum_{k=1,\dots,K}\mathbb{I}\{\kappa_k(s_{n+1})\ge\text{TL}\}\mathbb{I}\{\kappa_k(s_{n})<\text{TL}\}+\epsilon D,
\end{equation} 
where $D$ is the total time the machines where waiting for service in the period $\left(t_n, t_{n+1}\right]$ that is over the time limit $TL$ and is equal to
\begin{multline*}
D = \sum_{k=1,\dots,K}\mathbb{I}\left\{\kappa_k(s_{n+1}) > \text{TL}\right\}\min\left(t(s_{n+1})-t(s_n),\kappa_k(s_{n+1})-TL \right)+\\+\sum_{k=1,\dots,K}\mathbb{I}\left\{\kappa_k(s_{n})\ge 0 \right\}\mathbb{I}\left\{\kappa_k(s_{n+1})=-1\right\}\min\left(t(s_{n+1})-t(s_n), \kappa_k(s_{n})+t(s_{n+1})-t(s_n)-TL \right).
\end{multline*}

The first component of $D$ is the total waiting time over $TL$ during the transition period for all machines $k$ such that $\kappa_k(s_{n+1}) > TL$. The second component accounts for the event of type $e_t(s_{n+1})=4$ (a service engineer arrives at a machine), and adds the waiting time over $TL$ of the corresponding machine.

\section{Heuristic}\label{sec:heuristic}

In this section we describe the combined dispatching and relocation heuristic policy introduced in~\citep{pechina2019} for the problem described above that outperformed other heuristics studied in that paper. This heuristic consists of two parts. One is responsible for dispatching service engineers to broken machines, and the other part for relocation of idle service engineers between stations as well as allocation to base stations of service engineers that just became idle.\\

\noindent
\textbf{Dispatching}. First, we describe the dispatching part. One of the most common policies traditionally used in practice is the {\bf closest-first} policy, where the closest idle service unit is always dispatched to an incident, and the incidents in the queue are served in the first in, first out manner. However, in~\citep{pechina2019} it was empirically demonstrated that for our model always dispatching the closest idle service engineer whenever a failure occurs is not the best choice. It was shown that it can be beneficial to put a failed machine in a queue if there is a busy service engineer close enough to the breakdown. In that case the new breakdown is likely to be responded to quicker, and more idle service engineers are left in the system for future incidents.

The dispatching decision is made based on the notion of response time. For an idle service engineer $m$ response time $rt(k, m)$ to machine $k$ is the distance in time units from the current destination $l_m$ to the machine location $l_k$ plus the remaining distance to the current destination $d_m$. That is, $rt(k,m) = ||l_ml_k||_2 + d_m$. For a busy service engineer $m$ the response time is a random variable $RT(k,m) = ||l_ml_k||_2 + d_m + T^{rep}$, where $T^{rep}$ is an exponentially distributed repair time. When making a decision, the response time $RT(k,m)$ is estimated in such a way that the realization is likely to be less than the estimation. The estimation is done using an $\alpha^{\text{th}}$ percentile of the repair time distribution. Throughout the computational study, we take $\alpha = 80\%$. Thus, for a busy service engineer $m$ the response time to machine $k$ is estimated as $rt(k,m) = ||l_ml_k||_2 + d_m + T^{rep}_{80\%}$. This ensures the $80\%$ probability that the real response time is not larger than the estimation.

The dispatching policy works as follows. If the service engineer $m = \argmin_n rt(k, n)$, that minimizes response time, is idle, then that service engineer is dispatched immediately. If the service engineer $m$ is busy, then the call is placed in the queue. When a busy service engineer finishes a repair, he/she is dispatched to the closest machine from the queue, unless the queue is empty.\\

\noindent
\textbf{Relocation}. The relocation of service engineers is done using the DMEXCLP heuristic introduced in~\citep{jagtenberg2015} and adjusted for our model. This heuristic is used either when an idle service engineer is dispatched to a new breakdown, or when a service engineer finishes a repair and the queue is empty. In the first case, the heuristic considers relocating one of the rest idle stationary service engineers to another base station. In the second case, the heuristic allocates the service engineer that just finished the job to one of the base stations.

The DMEXCLP algorithm uses the notion of the \textit{expected covered demand}. To estimate the expected covered demand, we use the procedure introduced by Larson \cite{larson1974, larson1975}, and apply it to the model described in Section~\ref{sec:model}. Given the locations of the service engineers, the expected covered demand estimates the fraction of new requests that will be answered in time, and is approximated by
\begin{equation}
\frac{1}{|\mathcal{W}(s)|}\sum_{k\in \mathcal{\mathcal{W}(s)}}\sum_{i\in \mathcal{M}}P_i z_{ki}(s),
\label{eqn::expcovdem}
\end{equation}
where $\mathcal{W}(s)$ is the set of working machines in state $s$, $z_{ki}(s)$ is a binary variable indicating if the $i^\text{th}$ closest service engineer covers machine $k$ in state $s$, and $P_i$ is the probability that the first $i-1$ closest engineers are busy and the $i^\text{th}$ closest service engineer is idle. Note that machine $k$ is considered covered if there is at least one service engineer $m$, such that $rt(k, m) \leq TL$. The exact procedure for estimating $P_i$ is described in~\citep{pechina2019}.

According to the DMEXCLP algorithm, the action that maximizes the expected covered demand~\eqref{eqn::expcovdem} is always chosen, given that it satisfies certain restrictions. In the case when a service engineer becomes idle, that service engineer is allocated sequentially to each base station that is reachable within a certain time threshold $T_1$, then the resulting expected covered demand is computed using \eqref{eqn::expcovdem} for the obtained configuration, and the base station that leads to the best result is chosen. If there are no base stations within the given time threshold $T_1$, then all base stations are considered. In the case when a service engineer is dispatched to a new breakdown, we consider all pairs of base stations $\left(r_1, r_2\right)$ that are within the $T_2$ time units from each other, and with at least one service engineer at the base station $r_1$. The improvement in the expected covered demand is computed upon relocation of a service engineer from station $r_1$ to station $r_2$. If this improvement is larger than a given parameter $\Delta$, the relocation is made. If the gain in the expected covered demand is smaller than $\Delta$, or if there are no pairs of base stations within $T_2$ traveling time from each other, no relocation is made. In the computational study in Section \ref{sec:num:adp} we tune the parameters $T_1$, $T_2$ and $\Delta$ separately for each system using grid search.

\section{Approximate Dynamic Programming}\label{sec:adp}
The ADP approach combines the ideas of Markov decision theory and the heuristic approach. The goal is to approximate the value function as a linear combination of several so-called \textit{basis functions}, and choose actions using by substituting this approximation into the optimality equations. The idea is that the resulting policy will be close to optimal if the value function approximation is. 

Let us first formulate the optimality equations. Consider the process $\{s_n, \ n=0,1,\dots\}$ introduced in Section \ref{sec:model} with a discount factor $0<\gamma<1$. Let $V_{\pi}(s)$ denote the expected total discounted costs under policy $\pi$ when starting in state $s$:
\begin{equation*}
V_{\pi}(s) = \mathbb{E}\left[\sum_{n=0}^\infty \gamma^{t(s_n)}c(s_n, \pi(s_n), s_{n+1}) \mid s_0=s \right].
\end{equation*}
If policy $\pi^*$ is the optimal policy, then $V_{\pi^*}(s)$ satisfies the Bellman optimality equation
\begin{equation*}
V_{\pi^*}(s) = \min_{a \in \mathcal{A}(s)}\left\{\mathbb{E}_a\left[c(s, a, s') + \gamma^{t(s')-t(s)}V_{\pi^*}(s')\right]\right\}, \quad \forall s \in S,
\end{equation*}
where $s' = \Phi\left(s, a, \omega(s, a)\right)$ is the next state of the process when action $a$ is taken. We denote
\begin{equation*}
\pi^*(s) = \argmin_{a \in \mathcal{A}(s)}\left\{\mathbb{E}_a\left[c(s, a, s') + \gamma^{t(s')-t(s)}V_{\pi^*}(s')\right]\right\}, \quad \forall s \in S.
\end{equation*}
Note that there might be multiple minimizing actions in the above expression, leading to multiple optimal policies. In that case, $\pi^*$ refers to an arbitrary optimal policy.
In the remainder of this paper, we denote $V(s) := V_{\pi^*}(s)$ for ease of presentation.

To overcome the problem of a large (infinite in our case) state space, the ADP approach suggests to use an approximation $\hat{V}(s)$ of $V(s)$ that can be computed explicitly for any state $s$. We approximate $V$ by $\hat{V}$, a linear combination of several basis functions $\varphi_i(\cdot), \ i = 1,\dots,I$, i.e.,
\begin{equation*}
\hat{V}(\pmb{\alpha}, s) = \alpha_0 + \sum_{i=1}^I\alpha_i \varphi_i(s),
\end{equation*}
where $\pmb{\alpha} = (\alpha_0,\ldots,\alpha_I)$ is a vector of coefficients, and the approximate optimal policy is defined by
\begin{equation}\label{eq:adp:policy}
\hat{\pi}(\alpha, s) = \argmin_{a \in \mathcal{A}(s)}\left\{\mathbb{E}_a\left[c(s, a, s') + \gamma^{t(s')-t(s)}\hat{V}(\alpha,s')\right]\right\}, \quad \forall s \in S.
\end{equation}

As the state space of the process $\{s_n, \ n=0,1,\dots\}$ is countably infinite, the optimal action can not be computed in advance for each state. But if for each action $a \in \mathcal{A}(s)$ we could compute $\mathbb{E}_a\left[c(s, a, s') + \gamma^{t(s')-t(s)}V(s')\right]$ offline, then being in state $s$ we can find the optimal action $a$, as the action space is finite. In practice, the expected costs $\mathbb{E}_a\left[c(s, a, s') + \gamma^{t(s')-t(s)}V(s')\right]$ can be still hard to compute, even though the distribution of $s'$ depending on $s$ and $a$ is known (see Section \ref{sec:trans}), given that the state space is infinite. We estimate the value of $\mathbb{E}_a(\cdot)$ using Monte Carlo simulation, where the next state $s'$ is sampled $G$ times.
Then for each realization $s'_g,\ g=1,...,G$ we compute the future costs $c(s, a, s'_g) + \gamma^{t(s'_g)-t(s)}V(s'_g)$ and use the average of these costs as an approximation of the expected future costs. In our computational experiments we use $G=30$. 

The choice of basis functions is very important for the performance of the approach. First, they should be straightforward to compute for any state $s$. Second, they should jointly capture characteristics of the optimal value function, in order to obtain an accurate approximation. We discuss our choice of basis functions in Section~\ref{sec:adp:bf}. Given a set of basis functions, the approximation is finalized by tuning the vector of coefficients $\pmb{\alpha}$. We propose two metaheuristics to do this: a genetic algorithm and tabu search. We discuss these approaches in Section~\ref{sec:adp:tune}.

\subsection{Basis functions} \label{sec:adp:bf}

In this section we discuss the basis functions we consider for our model. Some of them describe the ability of the system to respond to future breakdowns (e.g., the number of uncovered machines and expected covered demand), while others capture future penalties for decisions made in the past (e.g., the number of unassigned calls and the number of unreachable calls).

{\bf Number of unreachable machines.}
Consider a machine for which a service engineer was already dispatched but he/she is still on the way and is not going to reach that machine in time. If its downtime did not yet exceed the time limit $TL$, we did not yet incur a penalty for this breakdown, but will in the future. The first basis function $\varphi_1(\cdot)$ represents the number of such machines:
$$\varphi_1(s) = \bigg|\Big\{k \in \mathcal{K} \mid 0 < \kappa_k(s) < \text{TL} \Big\} \cap \Big\{l_m(s), \ m \in \mathcal{M} \Big\} \bigg|.$$

{\bf Number of unassigned requests.}
Each unassigned repair request in the current state may result in future costs. First, this request may be responded to late and lead to penalty. Second, when a service engineer is going to be dispatched to the machine, the remaining coverage will decrease, which may in turn lead to costs for future simultaneous breakdowns. The second basis function counts the number of unassigned repair requests in state $s$:
$$\varphi_2(s) = \bigg|\Big \{k \in \mathcal{K}\setminus \left\{l_m(s), \ m \in \mathcal{M} \right\} \mid (\kappa_k(s) > 0) \ \lor \ (e_t(s) = 1)\land (e_l(s) = l_k)\Big \}\bigg|. $$

{\bf Number of missed unassigned requests.}
The missed requests are those with waiting time larger than the time limit $TL$. The missed unassigned requests are already incurring costs. Such requests also still require dispatching of a service engineer, which will lead to a decrease in coverage. At the same time as part of the costs is already incurred, such requests cannot be considered equal to other unassigned requests. The number of unassigned requests in state $s$ that already passed the time limit is defined as follows:
$$\varphi_3(s) = \bigg|\Big\{k \in \mathcal{K}\setminus \left\{l_m(s), \ m \in \mathcal{M} \right\} \mid \kappa_k(s) \ge \text{TL} \Big\} \bigg|. $$

{\bf Number of uncovered machines.}
A machine is considered \textit{covered} if there is at least one service engineer with the estimated response time less then $TL$ time unites. For each pair of machine $k$ and service engineer $m$ the response time $rt(k,m)$ is computed according to the procedure described in Section~\ref{sec:heuristic}. If a machine is not covered in state $s$ and a failure occurs, then this will result in costs, so the number of uncovered machines is an important metric. We consider only working machines as only those can break down:
$$\varphi_4(s) = \bigg|\Big\{k \in \mathcal{K} \mid (\kappa_k = 0) \land (l_k \neq e_l(s)) \land (rt(k, m) > \text{TL} \ \forall m \in \mathcal{M}) \Big\}\bigg|.$$

{\bf Expected covered demand.}
In Section~\ref{sec:heuristic}, we describe the relocation policy based on maximizing the expected covered demand that proved to perform well in~\cite{pechina2019}. Therefore, we include the expected covered demand in the set of basis functions:
$$\varphi_5(s) = \frac{1}{|W(s)|}\sum_{k \in \mathcal{W}(s)}\sum_{i\in \mathcal{M}}P_i z_{ki}(s).$$

{\bf Average response time.}
The last function is the average response time to the working machines. First, the response time $rt(k, m)$ is estimated for each pair of service engineer $m$ and machine $k$. Then for each demand node we choose the smallest response time over all service engineers and calculate the average over all machines. So, if $\mathcal{W}(s)$ is the set of working machines in state $s$, then
$$\varphi_6(s) = \frac{1}{|W(s)|}\sum_{k \in \mathcal{W}(s)}\min_{m\in \mathcal{M}}rt(k, m).$$

{\bf Future basis functions.}
All basis functions described above characterize the {\it current} state of the system. However, when we make a dispatching or relocation decision, we are not interested only in the state right after the decision is made, but also in the state upon arrival of the relocated repairman. Following~\cite{nasrollahzadeh2018}, for each basis function $\varphi_i(s)$ we introduce two basis functions $\varphi_i^{(1)}(s)$ and $\varphi_i^{(2)}(s)$ which characterize the state of the system after the arrival of the \textit{closest} and the \textit{second-closest} (respectively) traveling service engineer to their destinations.

Let $s^{(1)}$ denote the state of the system after the arrival of the closest traveling service engineer to his/her destination, and $s^{(2)}$ the state of the system after the arrival of the first two closest traveling service engineers to their destinations, assuming that the system is initially in state $s$ and no other events occur. The future basis functions are computed as $\varphi_i^{(1)}(s) = \varphi_i(s^{(1)})$ and $\varphi_i^{(2)}(s) = \varphi_i(s^{(2)})$. Note that if there are no traveling repairmen, then $\varphi_i^{(1)}(s) = \varphi_i^{(2)}(s) = \varphi_i(s)$.

\subsection{Tuning the ADP policy}\label{sec:adp:tune}

Standard methods for tuning the ADP policy aim at fitting the coefficients $\pmb{\alpha}$ so that $\hat{V}(\pmb{\alpha}, s) \approx V_{\pi^*}(s)$ for each state $s$, in hope that the obtained policy shows close to optimal performance. In this case, the coefficients $\pmb{\alpha}$ in~\eqref{eq:adp:policy} are chosen such that the distance between the optimal value function $\pmb{V}_{\pi^*}$ and $\pmb{\hat{V}_\alpha}$, the total discounted costs under the ADP policy with coefficients $\pmb{\alpha}$, is minimized:
\begin{equation}\label{eq:api}
\min_{\pmb{\alpha}}||\pmb{V}_{\pi^*} - \pmb{\hat{V}_\alpha}||_p.
\end{equation}
One of such algorithms is called \textit{approximate policy iteration}, where the current policy value function is evaluated in a simulation, and the coefficients are iteratively updated using linear regression, that is, the value $p=2$ is used in~\eqref{eq:api}. As noted in~\cite{maxwell2013}, despite some of the advantages of approximate policy iteration (e.g., ease of understanding and implementation), it also has drawbacks that might lead to low quality solutions in terms of actual performance. We implemented approximate policy iteration, but observed consistently poor performance despite convergence to low values of the mean squared error in the linear regression. Note that the choice of an action in a given state under the ADP policy~\eqref{eq:adp:policy} depends on the relative difference in value function for different states, rather than on the actual values. As our goal is to obtain a high performance policy rather than to fit the value function in~\eqref{eq:api}, we resort to heuristics that tune the coefficients $\pmb{\alpha}$ based on the corresponding performance in a simulation.

\subsubsection{Genetic algorithm} \label{sec:adp:ga}
The first approach we propose is a genetic algorithm. Genetic algorithms draw inspiration in evolution and natural selection, and are widely used for optimization problems \cite{yu2010}. The key concepts of genetic algorithms are
\begin{itemize}
\item \textit{Population}: each individual/solution is a part of a pool;
\item \textit{Selection}: the fittest individuals survive;
\item \textit{Crossover}: the fit individuals reproduce, propagating their fit genes;
\item \textit{Mutation}: sometimes new characteristics appear by accident.
\end{itemize}

In our case an individual is a vector of coefficients $\pmb{\alpha}$, and a population is a set of vectors $\pmb{\alpha}^{(n)}, \ n = 1, \dots, N^{GA}$, where $N^{GA}$ is the size of the population and one of the parameters of the algorithm. In each iteration of the algorithm, the population is updated using \textit{mutation}, \textit{crossover} and \textit{selection}. The population is initialized by adding random vectors to the same vector $\alpha^{(0)}$. In our experiments we set all elements of vector $\alpha^{(0)}$ equal to $1$ except those corresponding to the expected covered demand, for which the coefficient is set to $-1$. To produce an individual of an initial population, we then add a random value to each element of $\alpha^{(0)}$ drawn from a uniform distribution $\mathcal{U}(-1, 1)$.

During each iteration of the algorithm, a new population is constructed as follows. First, the {\it crossover} operator is used, where $N^{GA}$ pairs of candidate solutions are randomly selected from the current population and the average of each pair is chosen as a candidate for the next population. Second, the {\it mutation} operator is applied to each candidate solution from the current population. Given a solution $\pmb{\alpha}^{(n)}$, the {\it mutation} operator adds a vector, the $i$th component of which is normally distributed with mean $0$ and standard deviation $A^{GA}|\pmb{\alpha}^{(n)}_i|$. Here $A^{GA}$ is the parameter of the algorithm, and is referred to as \textit{mutation amplitude}. This way, together with the current population, after mutation and crossover is done, we obtain $3N$ candidate solutions for the next population.

Finally, the selection is done to produce the next population, where $N^{GA}$ solutions are chosen out of the $3N^{GA}$ candidates. To measure fitness of a candidate solution, we run simulation with the corresponding ADP policy. Simulation starts from the same state of each candidate, and the fitness of a given candidate solution is measured as the costs per machine failure observed in a simulation (see the cost structure in Section~\ref{sec:costs}). Note that the larger the time horizon of simulation, the better is an estimation of the solution fitness, but the larger the computational time of the algorithm. In our experiments we set the time horizon equal to $500/(\lambda K)$.

\begin{figure}[h]
\centering
\includegraphics[width = \textwidth]{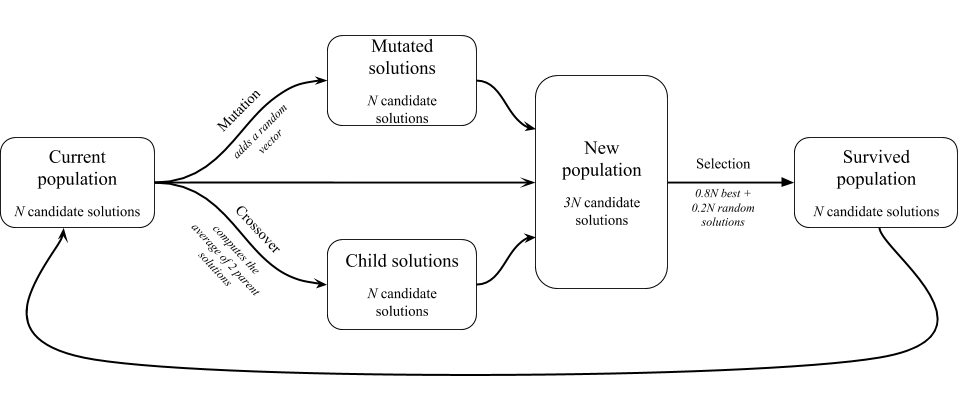}
\caption{Genetic algorithm.}
\label{fig::gen_alg}
\end{figure}

All $3N$ solutions of the new population are evaluated by means of simulation: the system is simulated under corresponding policies from the same initial state and the fraction of calls answered in time is observed. Then only the $q^{GA}N^{GA}$ best performing and $(1-q^{GA})N$ randomly chosen candidate solutions survive and constitute the next population. Here the fraction of fittest candidate solutions $q^{GA}$ is the parameter of he algorithm. The algorithm stops after a certain number of iterations. The scheme of the algorithm is depicted in Figure \ref{fig::gen_alg}.

\subsubsection{Tabu search} \label{sec:adp:ts}

Tabu search is a high-level metaheuristic technique that guides other local search heuristic methods and is constructed in such a way that allows to escape the local optima~\citep{glover1998}. The main idea is to allow moves to worse solutions and prevent cycling back to the local optima with the help of the so called tabu lists. Here we introduce the main concepts of the tabu search metaheuristic:

\begin{itemize}
\item \textit{Incumbent solution}: the current solution representing the current state of the algorithm. Tabu search performs a walk in the search space through a sequence of incumbent solutions. The best found solution is the outcome of the algorithm.
\item \textit{Move}: a procedure of obtaining a new feasible solution from the incumbent solution.
\item \textit{Neighborhood}: determined by the move and represents all the feasible solutions that can be reached by moving from the incumbent one.
\item \textit{Tabu list}: a list of restrictions that impose limitations based on certain attributes of the recently performed moves. The tabu list has a limited size, it is based on a short term memory and prevents the search to return to the recently visited solutions.
\item \textit{Tenure} is a number of iterations of the algorithm a tabu restriction stays in the tabu list.
\end{itemize}

The algorithm has two phases. First, it performs a certain number of iterations in the diversification phase, followed by a certain number of iterations in the intensification phase that starts from the best found solution of the diversification phase. In the diversification phase, the focus is on exploring the neighborhood and escaping the local optima. In each iteration, the algorithm moves to the best found candidate solution from the neighborhood, even if it is worse than the current incumbent solution. In the intensification phase, the goal is to find the local optimum around the best found solution of the diversification phase. The movement to a new incumbent solution is allowed only if it is better than the current one.

The {\it move} operator works as follows. Given an incumbent solution $\pmb{\alpha}$, the {\it move} operator adds a random vector to it, $i$th component of which is normally distributed with mean $0$ and standard deviation $A^{TS}|\pmb{\alpha}_i|$, unless $i$ is in the tabu list. The tabu list contains the indices of the coefficients in the incumbent solution that are not allowed to change, as well as the remaining number of iterations those indices will stay in the tabu list. The size of the tabu list $TLS$ is one of the parameters of the algorithm. The tabu list is initialized as an empty set. In the end of each iteration, the remaining number of iterations is decreased by one. If it becomes zero, a new set of $TLS$ indices is chosen as tabu with the remaining number of iterations equal to tenure $TLT$, another parameter of the algorithm. In the diversification phase, the tabu indices are chosen as non-tabu indices that changed the most in the last iteration, and in the intensification phase those that changed the least.

We initialize the algorithm with the incumbent solution $\alpha^{(0)}$ defined the same way as in the genetic algorithm. At every iteration $N^{TS}$ candidates are obtained from the incumbent solution using the {\it move} operator. The candidate solutions are estimated first by the surrogate fitness function, and then a subset of $f^{TS}N^{TS}$ good candidate solutions is estimated by the primary fitness function. Both fitness functions are measured as costs per machine failure observed in a simulation. The only difference is that the surrogate fitness is computed over a shorter time horizon to quickly identify potentially good candidates. Those are further evaluated in a larger simulation. Using surrogate fitness allows us to better explore the neighborhood at a lower computational cost. In our experiments, we use the time horizon of $30/(\lambda K)$ time unites for the surrogate fitness, and $500/(\lambda K)$ for the primary fitness. If more than $f^{TS}N^{TS}$ candidates result in zero surrogate fitness (may happen when the quality of the incumbent solution improves), the corresponding time horizon is doubled.

At the end of each iteration, the algorithm updates the incumbent solution and the tabu list. If in the intensification phase the incumbent solution does not change, then the {\it move} amplitude is decreased by $10\%$ before the next iteration. The algorithm stops after a certain number of iterations in the diversification and intensification phases. If necessary, the cycle can repeat starting from the incumbent solution. In the numerical experiments in Section~\ref{sec:num} we perform only one such cycle.

\section{Numerical results}\label{sec:num}

In this section, we present the setup and the results of our numerical experiments. In Section~\ref{sec:num:setup} we describe the types of networks we use together with the parameters defining the important properties of the networks. Next, in Section~\ref{sec:num:par} we study the parameters of the two ADP tuning algorithms presented in Section~\ref{sec:adp:tune}. Finally, in Section~\ref{sec:num:adp} we use simulation to compare the ADP policies obtained by both tuning algorithms against the heuristic described in Section~\ref{sec:heuristic} and the closest-first policy.

\subsection{Setup of the numerical experiments}\label{sec:num:setup}

In our computational experiments, the networks are generated randomly as follows. First, we define a square on a Euclidean plane. Then we randomly generate the coordinates of the base stations within that square. After that, the coordinates of the machines are sampled at random such that each machine is within the distance of $TL$ units from at least one station. We only consider the networks where each base station covers at least one machine. After the locations of the machines and the base stations are determined, the service engineers are allocated to the base stations such that the expected covered demand is maximized. This is done by solving an integer program~\citep{pechina2019}.

\begin{figure}[t]
\centering
        \begin{subfigure}[b]{0.5\textwidth}
                \includegraphics[width=\linewidth]{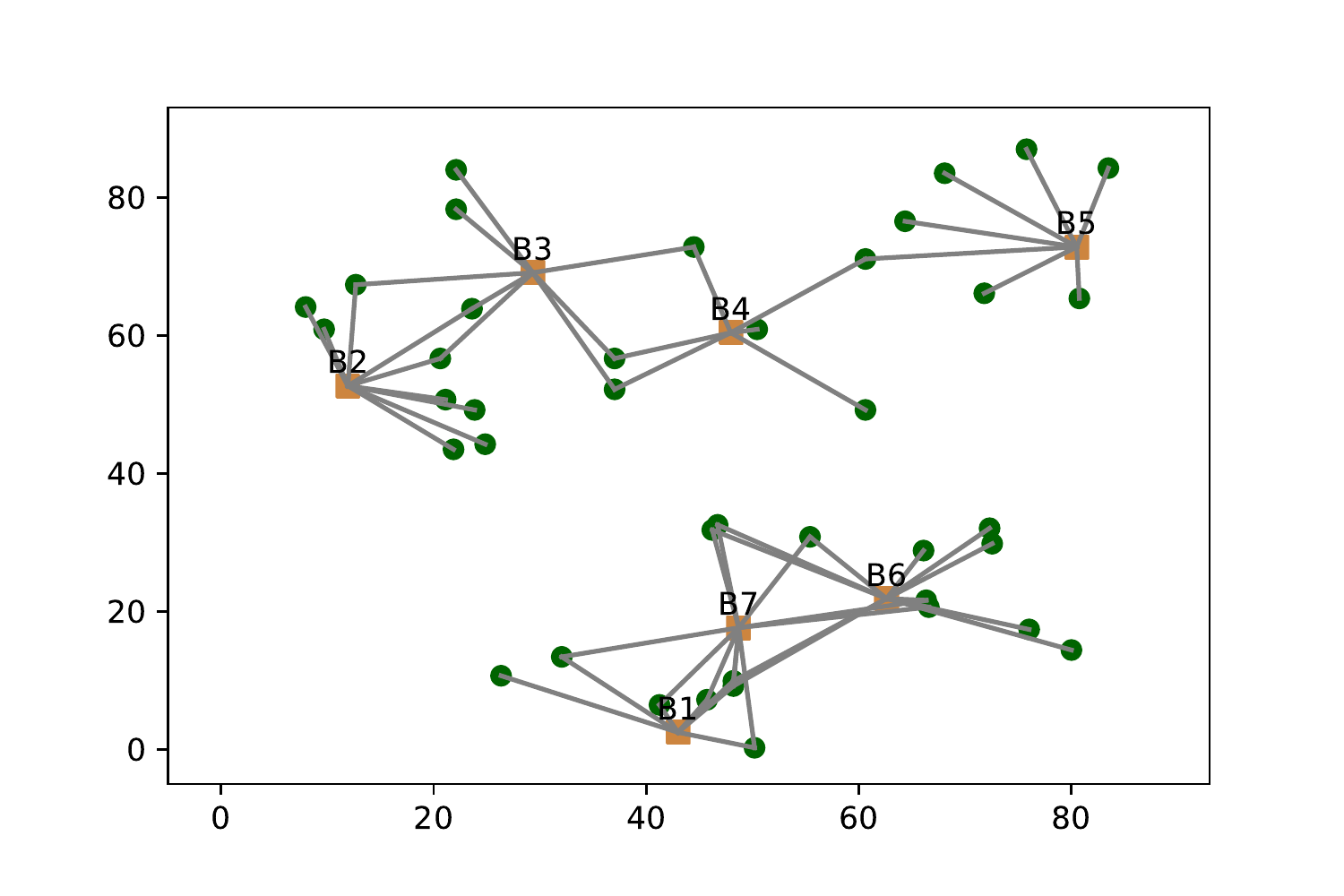}
                \caption{$d = 0.3$}
                \label{fig::density:0.3}
        \end{subfigure}%
        \begin{subfigure}[b]{0.5\textwidth}
                \includegraphics[width=\linewidth]{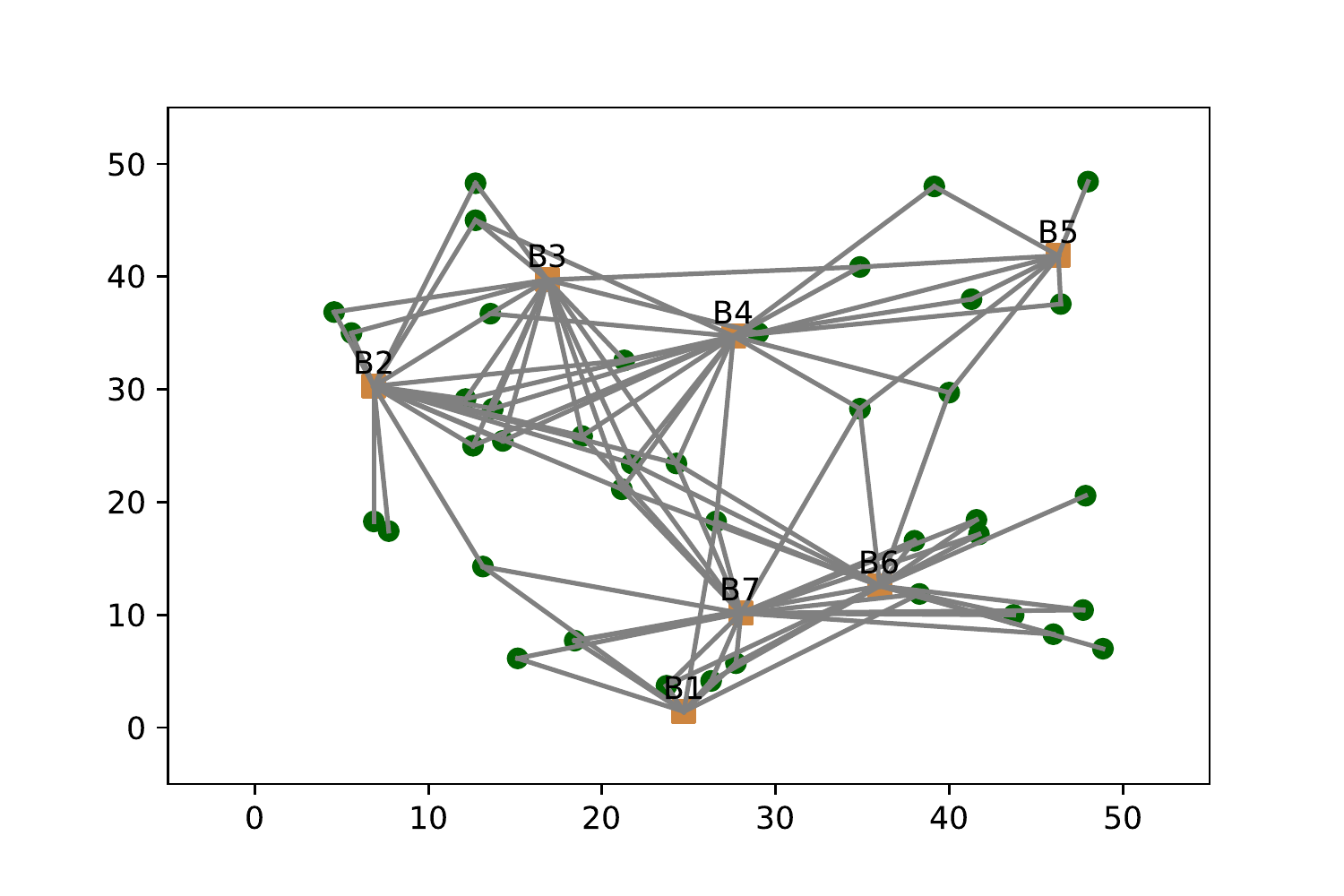}
                \caption{$d = 2$}
                \label{fig::density:2}
        \end{subfigure}%
        \caption{Networks of different density $d$, given $TL=20$, $K=40$, $R=10$.}\label{fig::density}
\end{figure}

The size of the square is determined by the combination of the time limit $TL$ and the density parameter $d$. The higher the density for a given value of $TL$, the smaller the square, meaning that the map is more dense with each node covered by more stations on average. An example of how $d$ affects the network structure can be seen in Figure~\ref{fig::density}. The two maps are randomly generated for two values of $d$ and the same values of $TL$, $R$ and $K$. The base stations are connected with an edge to the machines they cover. The map with low value of $d$ is more sparse. In sparse networks, the machines are covered by less stations, and the distances are larger compared to $TL$. On the one hand, if the distances are large, relocation might be not desirable as it takes more time. On the other hand in sparse networks the same station may cover multiple machines if there are more machines than stations. In that case, relocation from a full station to an empty one upon a incident may be beneficial for the system performance. The map with high value of $d$ has smaller distances and more connecting edges. If the map is dense, then the machines are located closer to each other, meaning that a busy service engineer may potentially respond quicker to a nearby incident than the closest idle one. In that case, the optimal policy might be putting the new incident in the queue instead of dispatching an engineer immediately.

The other two parameters affecting the policy performance are the breakdown rate $\lambda$ and the service rate $\mu$. These two parameters together define the load of the system; that is, the fraction of time the service engineers are busy responding to incidents. In our experiments, we fix the parameter $\lambda$ and control the load via $\mu$. The larger $\mu$ the faster the service engineers fix the breakdowns, and hence, the load is lower.

\subsection{Parameters of the ADP tuning algorithms}\label{sec:num:par}

\begin{figure}[t]
\centering
\includegraphics[width=0.7\linewidth]{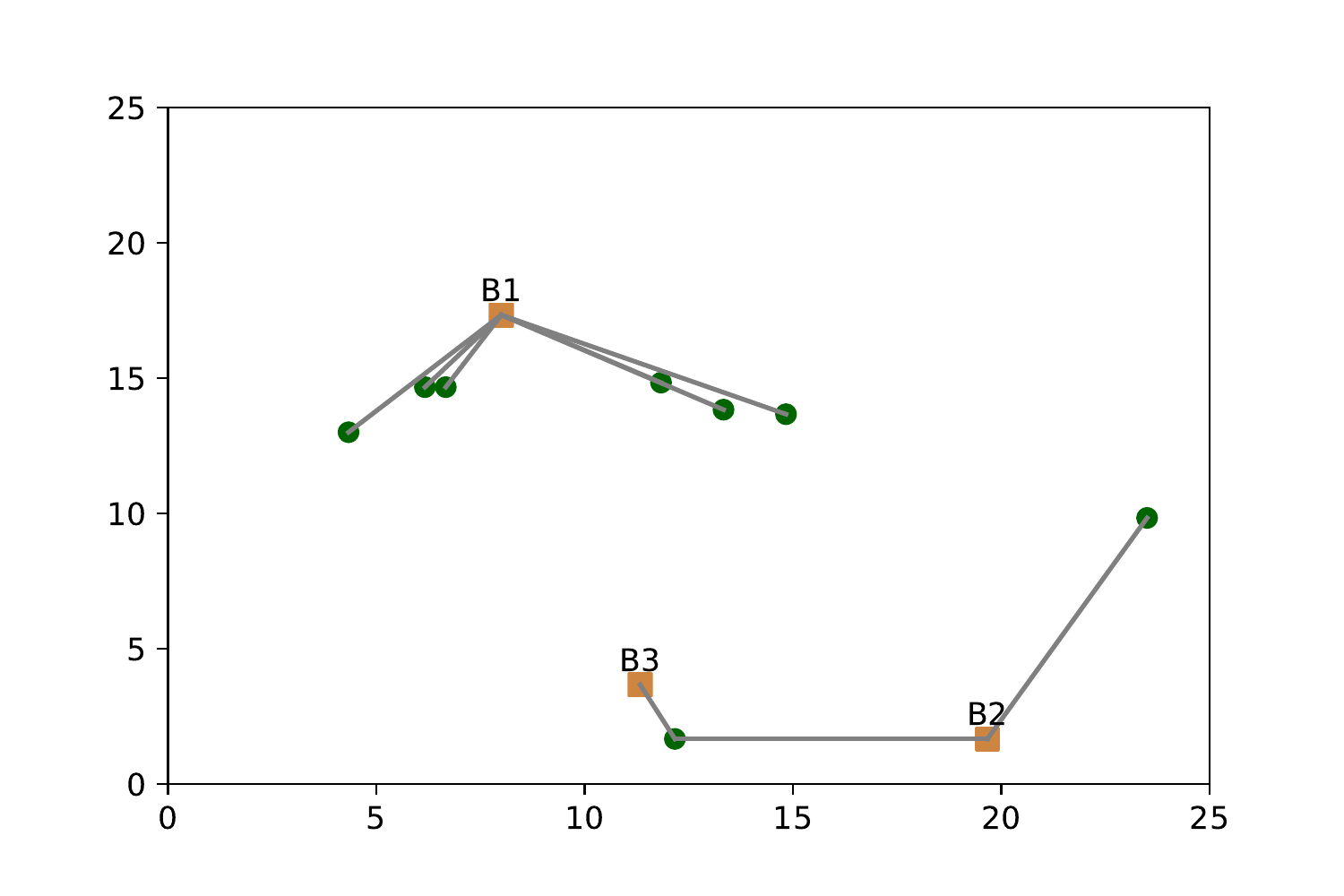}
\caption{Network used to test the parameters of the ADP tuning algorithms.}
\label{fig:tune_par:network}
\end{figure}

In this section, we discuss the important parameters of the two ADP tuning algorithms, genetic algorithm and tabu search, introduced in Section~\ref{sec:adp:tune}. There are two criteria that affect the choice of the values for each of the parameters. The first one is the convergence rate, that is, how quickly in terms of the number of iterations the algorithm finds good quality solutions. The second criterion is the computational time, that is, how much time the algorithm spends per iteration. The right balance should be found so that the good solution are found within a reasonable amount of time. In all experiments in the rest of this section we use the randomly generated network shown on Figure~\ref{fig:tune_par:network} with $K=8$, $R=3$, $M=5$, $TL=10$, $\lambda = 0.01$ and $\mu = 0.1$. The fitness of each solution is measured in terms of cost per incident incurred in a simulation.
\newline

\begin{figure}[t]
\centering
\includegraphics[width=0.7\linewidth]{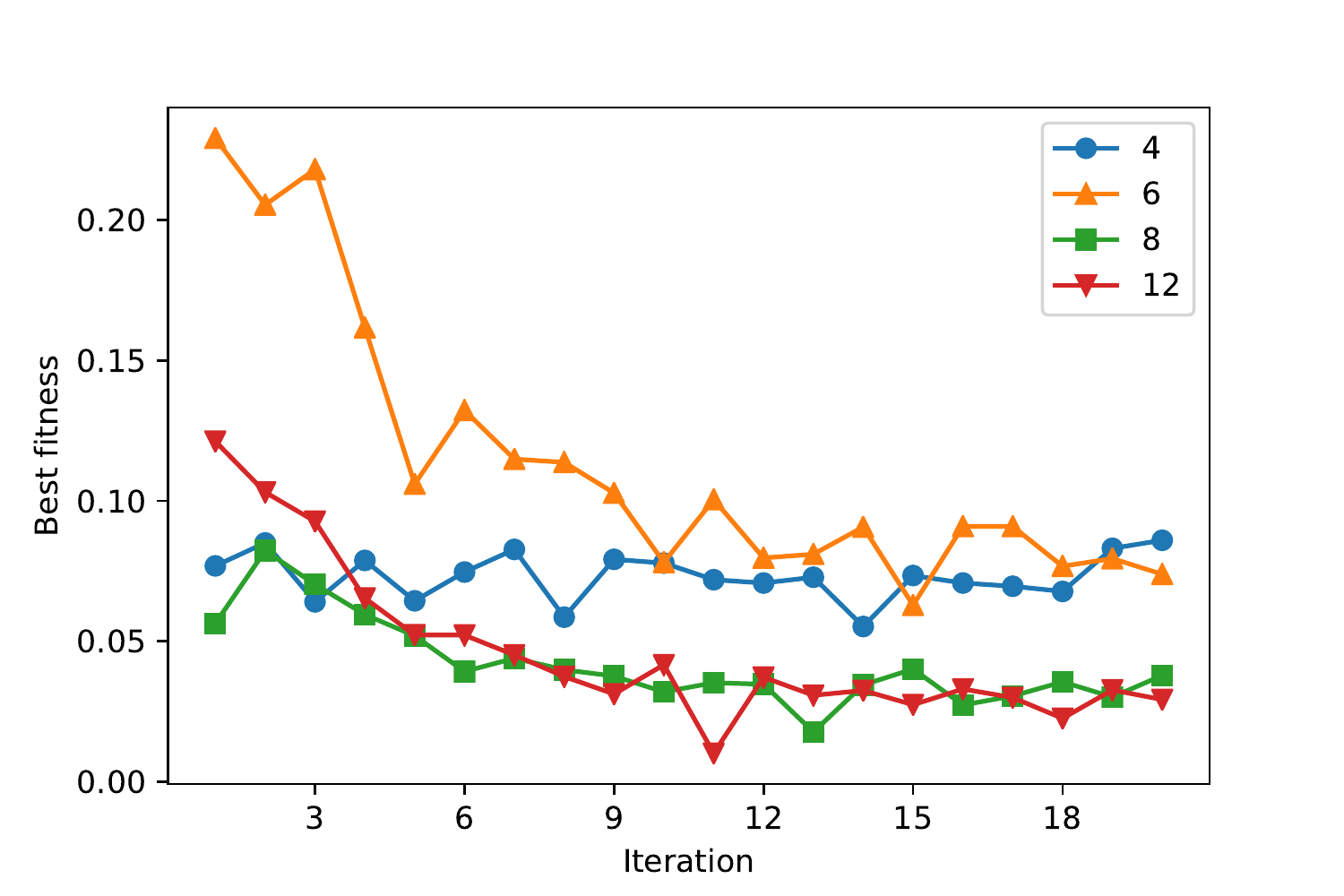}
\caption{Genetic algorithm convergence for different sets of basis functions. $N^{GA}=50$, $A^{GA}=1$, $q^{GA}=0.8$.}
\label{fig:tune_par:bf}
\end{figure}

\noindent
\textbf{Number of basis functions.}\\
First, we determine if there are basis functions that do not contribute to the quality of the solutions found. The more basis functions are included, the more complex the problem becomes, as each extra basis function adds an extra dimension to the search space. This leads to large computational times. In our experiments, we discovered that when an ADP tuning algorithm is ran with the six main basis functions $\varphi_1(\cdot),\ldots,\varphi_6(\cdot)$, the coefficients corresponding to the functions $\varphi_4(\cdot)$ and $\varphi_6(\cdot)$ are driven to zero, meaning that these two basis functions have insignificant effect on the quality of the ADP policy. Since the number of uncovered machines and the average response time seam to bare redundant information about the state of the system, we, therefore, consider omitting the corresponding main and future basis functions. Figure~\ref{fig:tune_par:bf} demonstrates the convergence of the genetic algorithm depending on the number of basis functions used. Here, 6 corresponds to the six main basis functions, 4 to the main basis functions excluding $\varphi_4(\cdot)$ and $\varphi_6(\cdot)$, 8 to the previous four plus the corresponding one-step-ahead future basis functions, 12 to the previous eight plus the corresponding two-steps-ahead future basis functions. The fitness of the best solution in a population is plotted against the number of iteration.

Note that the initial populations are different for different sets of basis functions, and the solutions are sampled at random. So, the quality of the initial population does not depend on the set of basis functions used. It can be seen from Figure~\ref{fig:tune_par:bf} that using six main basis functions instead of four does not contribute to the quality of the obtained policies. Adding the one-step-ahead future basis functions, however, boosts the performance of the policies. Including the two-steps-ahead future basis functions further improves the quality of the best found policy, although the effect is marginal. A similar effect is observed when running the tabu search algorithm. In Section~\ref{sec:num:adp} we use 12 basis functions with both algorithms when tuning the ADP policy.
\newline

\begin{figure}[t]
\centering
        \begin{subfigure}[b]{0.5\textwidth}
                \includegraphics[width=\linewidth]{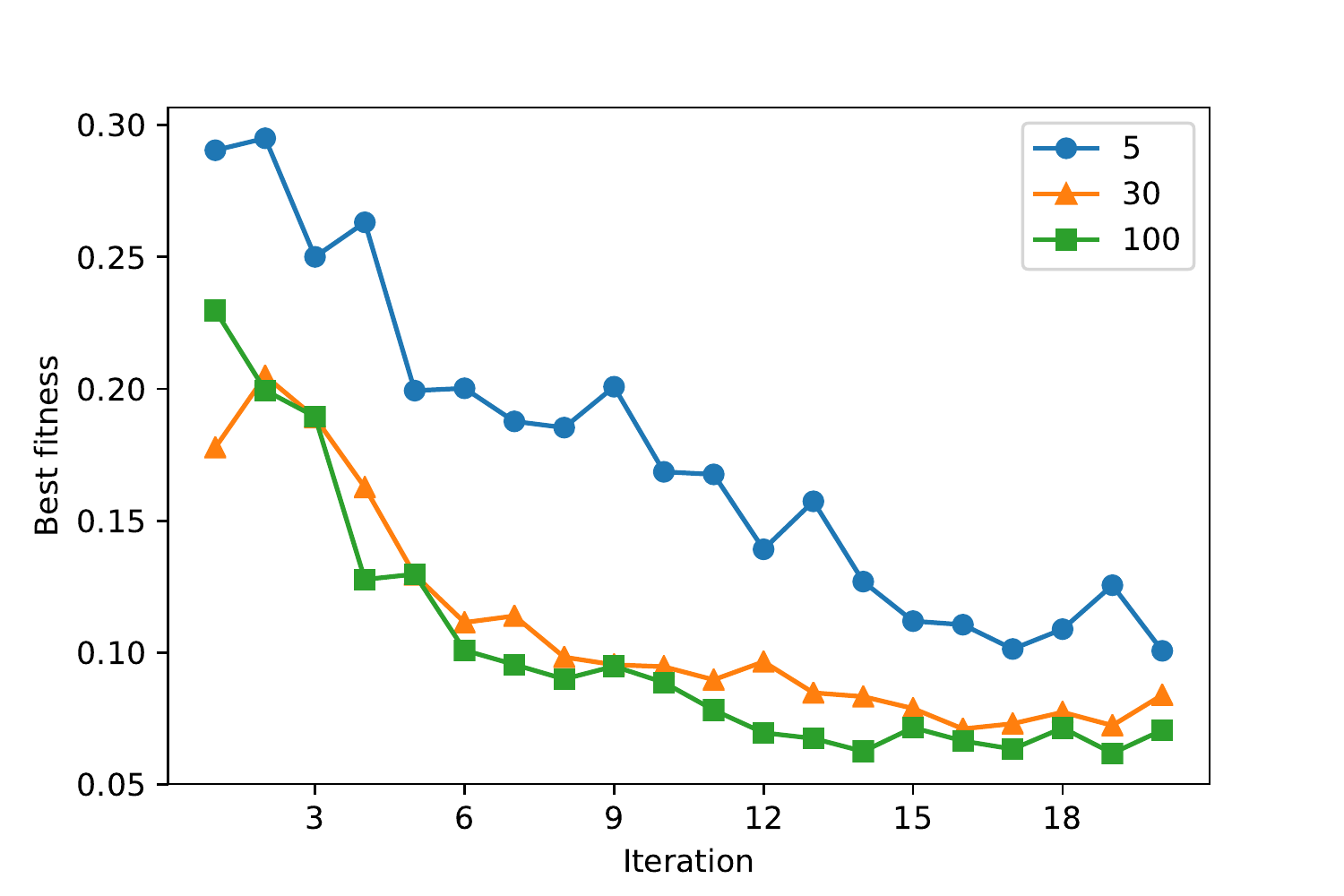}
                \caption{$\pmb{N^{GA}}$ ($A^{GA}=1$, $q^{GA}=0.8$)}
                \label{fig::density:0.3}
        \end{subfigure}%
        \centering
        \begin{subfigure}[b]{0.5\textwidth}
                \includegraphics[width=\linewidth]{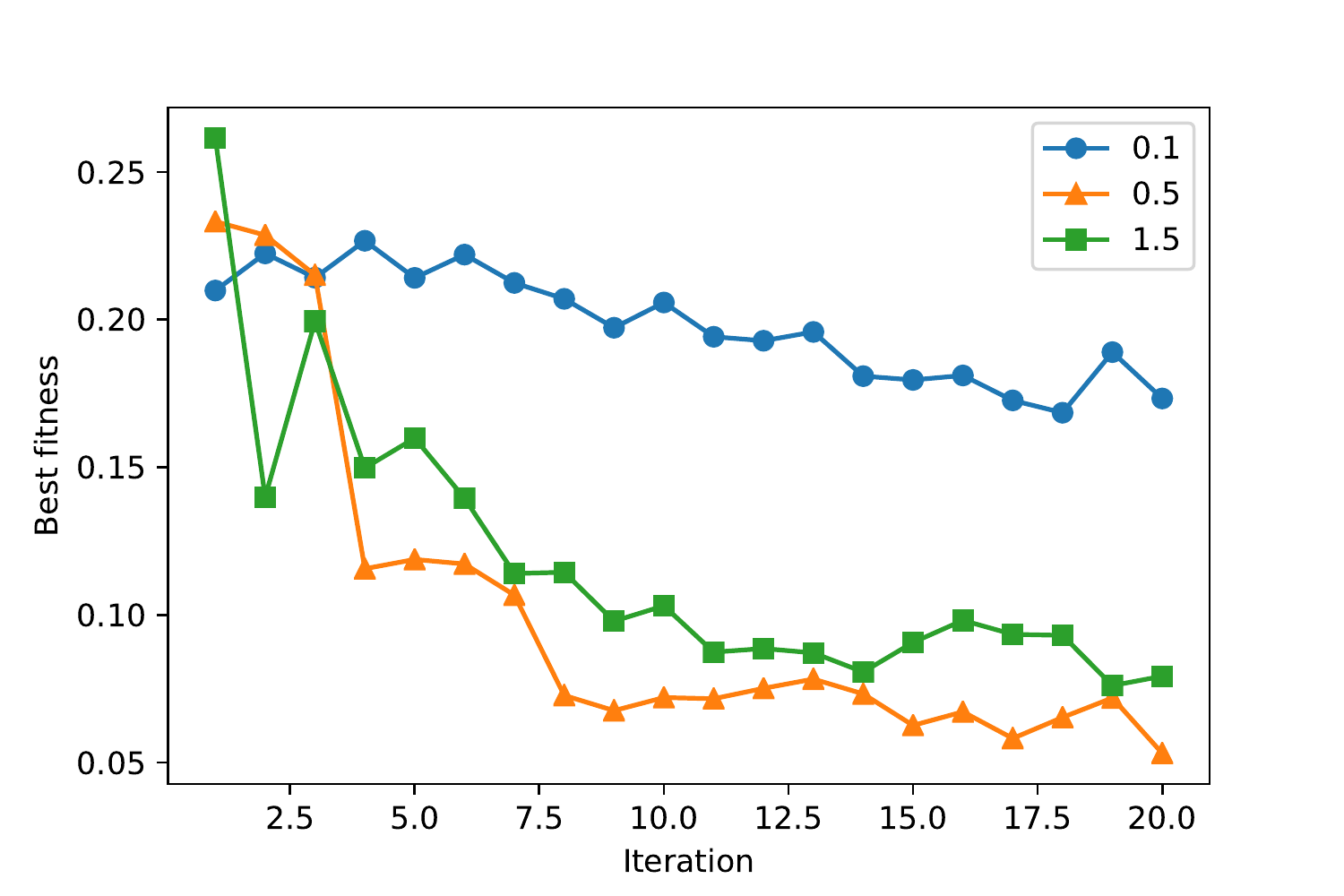}
                \caption{$\pmb{A^{GA}}$ ($N^{GA}=50$, $q^{GA}=0.8$)}
                \label{fig::density:2}
        \end{subfigure}
        \centering
        \begin{subfigure}[b]{0.5\textwidth}
                \includegraphics[width=\linewidth]{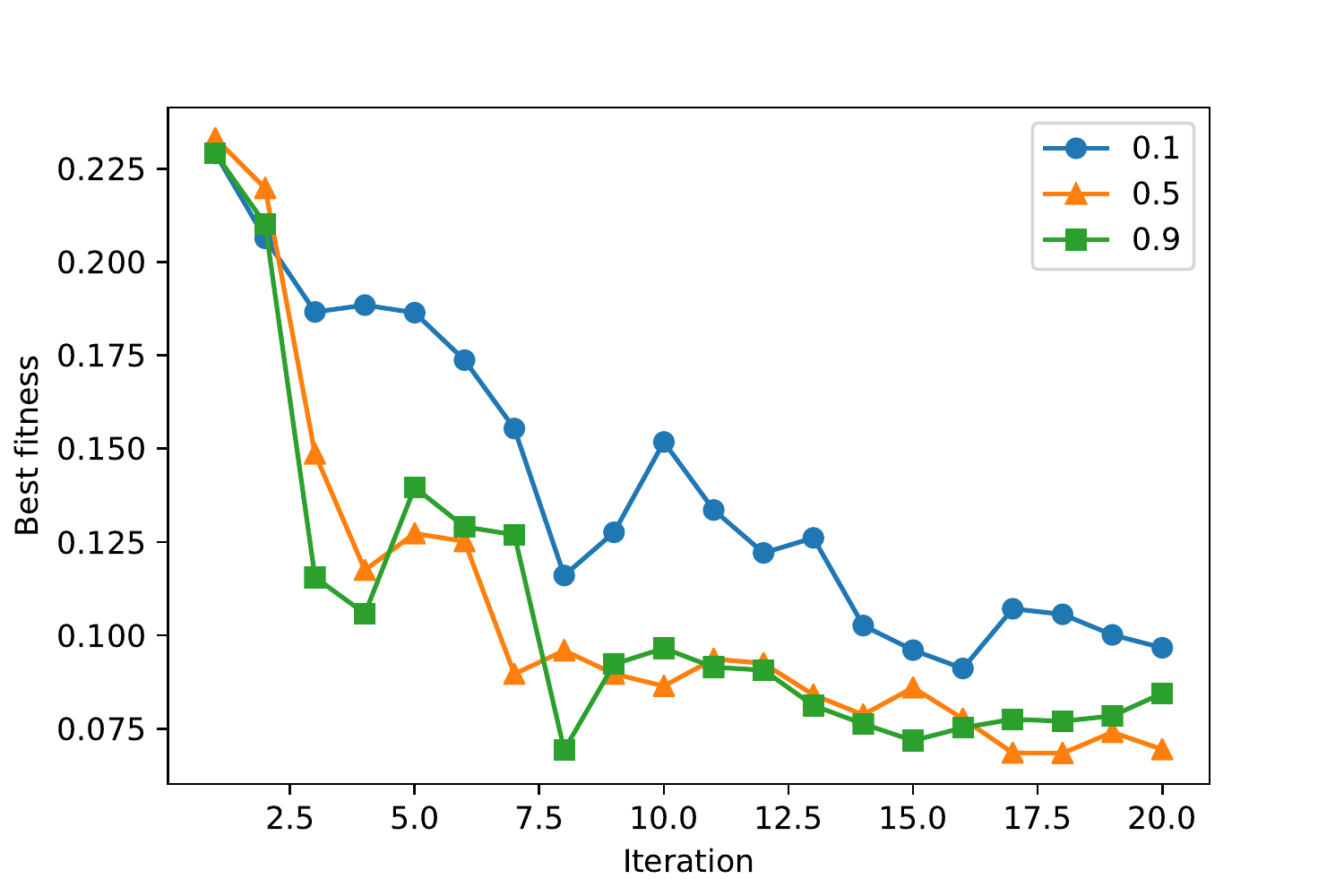}
                \caption{$\pmb{q^{GA}}$ ($N^{GA}=50$, $A^{GA}=1$)}
                \label{fig::density:2}
        \end{subfigure}%
        \caption{Genetic algorithm convergence with 6 main basis functions.}
        \label{fig:tune_par:ga}
\end{figure}

\noindent
\textbf{Parameters of genetic algorithm $\pmb{N^{GA}}$, $\pmb{A^{GA}}$, $\pmb{q^{GA}}$.}\\
Figure~\ref{fig:tune_par:ga} presents the effects of the genetic algorithm parameters on the fitness convergence, where the best population fitness is plotted per iteration. It can be seen that increasing the population size $N^{GA}$ has a positive effect. It, however, directly affects the computation time, as in each iteration a simulation is run for every solution in the population. Hence, the choice should be made depending on the computational resources to ensure that the policy can be obtained within a reasonable amount of time. In Section~\ref{sec:num:adp} we use $N^{GA}$ between 50 and 100 depending on the system. The mutation amplitude $A^{GA}$ and the fraction of fittest $q^{GA}$ do not affect computation time, but do influence the convergence of the algorithm. Both parameters should not be too high or too low. When tuning the ADP policy with genetic algorithm in Section~\ref{sec:num:adp}, we set $A^{GA}=0.7$ and $q^{GA}=0.8$.
\newline

\begin{figure}[t]
\centering
        \centering
        \begin{subfigure}[b]{0.49\textwidth}
                \includegraphics[width=\linewidth]{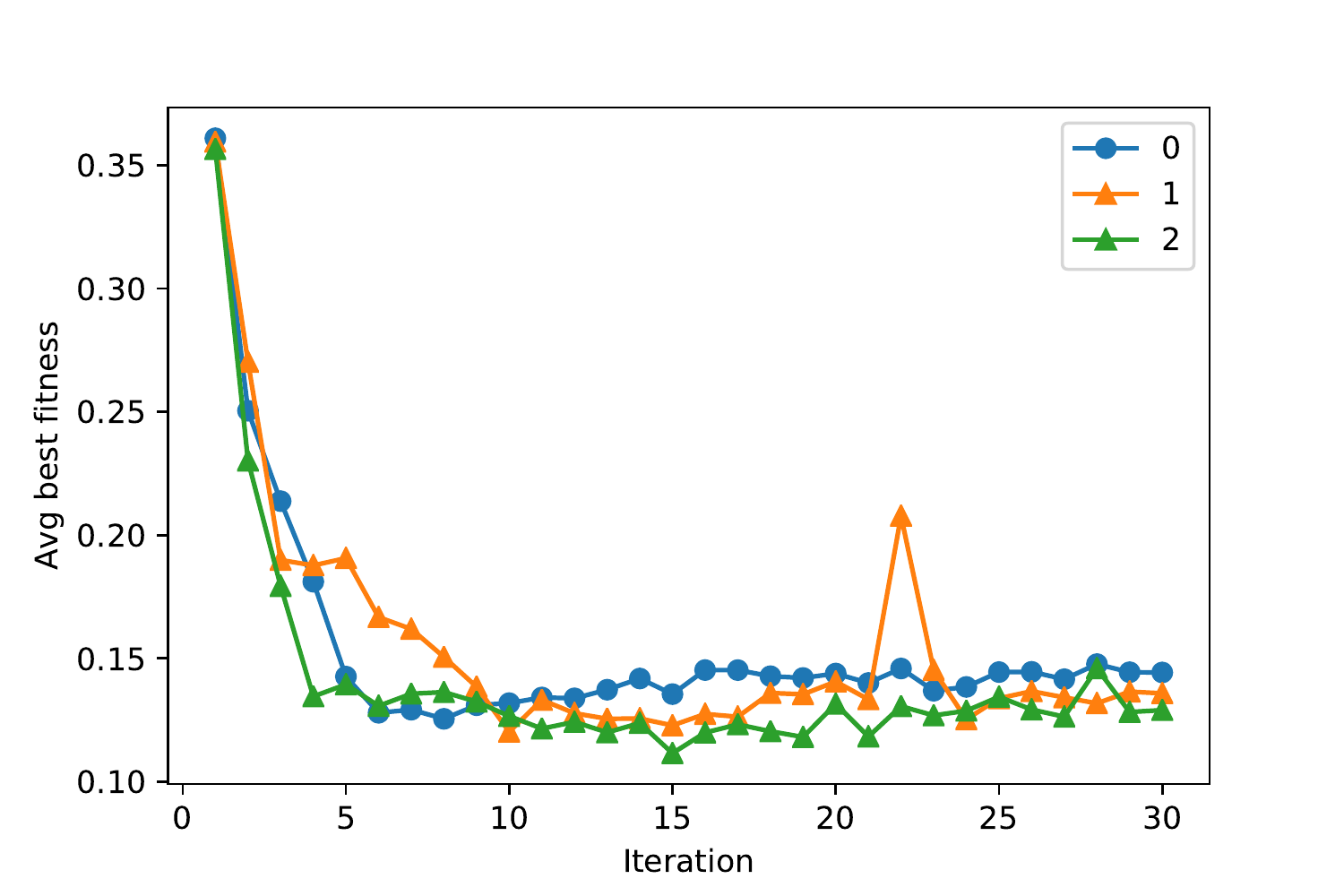}
                \caption{$\pmb{TLS}$ \\($N^{TS}=400$, $A^{TS}=1$, $f^{TS} = 0.3$, $TLT = 1$)}
                \label{fig::density:2}
        \end{subfigure}
        \centering
        \begin{subfigure}[b]{0.49\textwidth}
                \includegraphics[width=\linewidth]{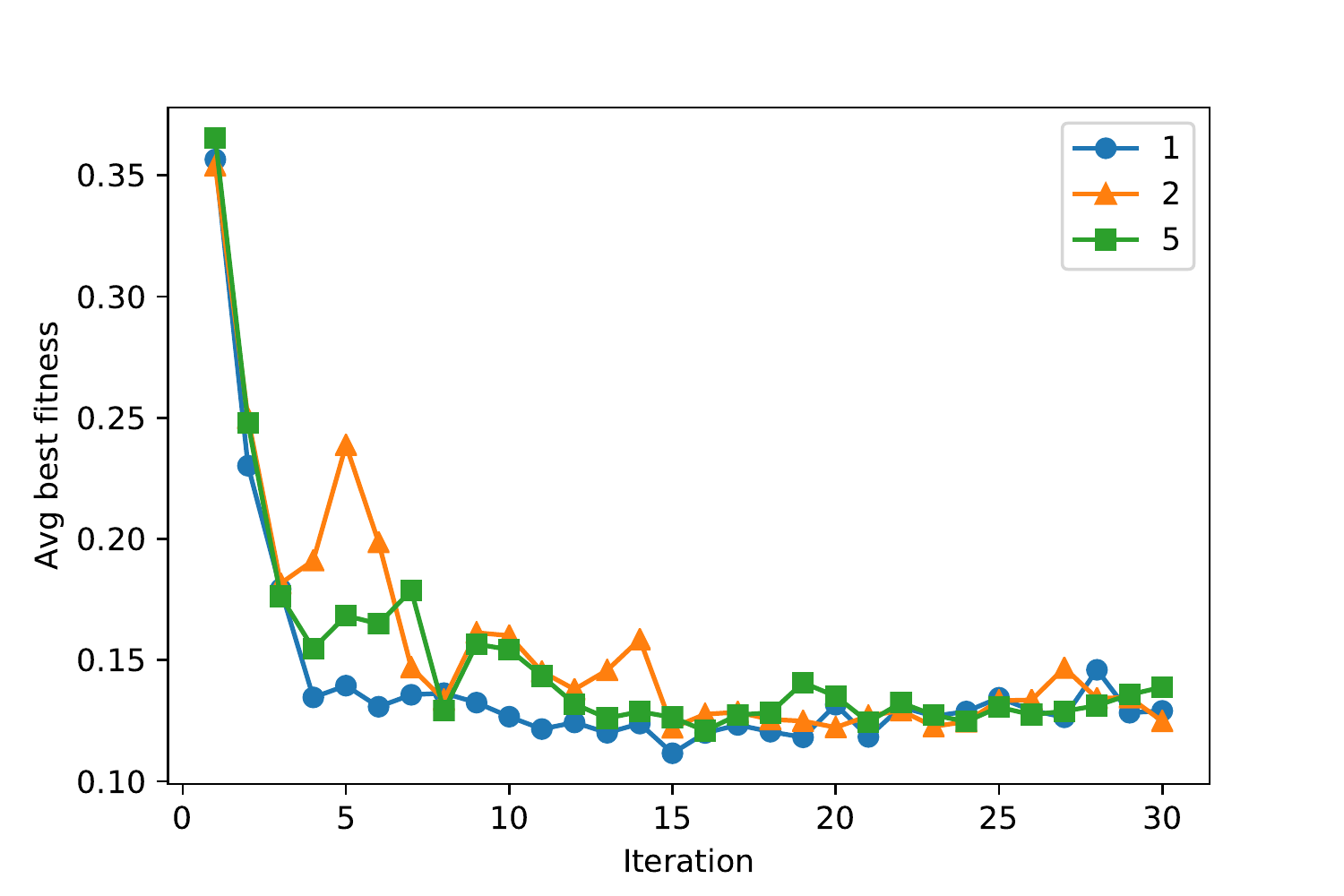}
                \caption{$\pmb{TLT}$ \\($N^{TS}=400$, $A^{TS}=1$, $f^{TS} = 0.3$, $TLS = 2$)}
                \label{fig::density:2}
        \end{subfigure}
        \caption{Tabu search convergence with 6 main basis functions.}
        \label{fig:tune_par:ts}
\end{figure}

\noindent
\textbf{Parameters of tabu search $\pmb{f^{TS}}$, $\pmb{TLS}$, $\pmb{TLT}$.}\\
The parameters $N^{TS}$ and $A^{TS}$ of the tabu search algorithm demonstrate similar influence as the $N^{GA}$ and $A^{GA}$ parameters of the genetic algorithm. Therefore, we do not discuss them here. In our computational experiments in Section~\ref{sec:num:adp}, we use the values between 200 and 400 for $N^{TS}$, and we set $A^{TS}=0.7$. The first tabu search specific parameter we consider is $f^{TS}$. It is intuitive that the larger its value, the better the convergence in fitness, as for larger values of $f^{TS}$ a bigger part of the neighborhood is explored in each iteration. It is also intuitive that this leads to larger computation time, as in each iteration more solutions from the neighborhood are evaluated with the primary fitness function. When fitting ADP policies in Section~\ref{sec:num:adp}, we set $f^{TS} = 0.3$. Figure~\ref{fig:tune_par:ts} shows how the parameters $TLS$ and $TLT$ affect the algorithm convergence, where the average primary fitness of the best candidates in the neighborhood is plotted per iteration of a diversification phase. These parameters have no effect on computation time, but do influence the convergence of the algorithm and should be chosen carefully. In the example shown in Figure~\ref{fig:tune_par:ts} increasing $TLS$ helped to find a better neighborhood, but at a cost of making more iterations. Increasing $TLT$, however, did not help converging to a good neighborhood. Note that the right choice of these parameters depends on the number of basis functions used. In Section~\ref{sec:num:adp}, we set $TLS=3$ and $TLT=2$ in combination with twelve basis functions.

\subsection{ADP performance}\label{sec:num:adp}

In this section, we compare the ADP policy against the heuristic policy (Heuristic) described in Section~\ref{sec:heuristic} and the closest-first (CF) policy that always dispatches the closest engineer and does not perform relocations. We consider various types of systems by changing the parameters $\mu$, $d$ and $TL$. We fix $\lambda$ at 0.01 and the size of the networks with $K=40$, $R=7$ and $M=10$, which is realistic for many real life maintenance networks. For each combination of the parameters 10 random networks are generated as described in Section~\ref{sec:num:setup}. For each of the networks the ADP policy is tuned with both genetic algorithm (GA) and tabu search (TS). The two ADP policies are then compared against the Heuristic and the CF policies using simulation over the time horizon of $1000/\lambda = 100000$ time units. The performance of the policies is measured with the cost per incident (computed according to~\eqref{eq:costs} with $\epsilon=0.01$ and referred to as Cost), fraction of late arrivals (FLAR), and average response time (ART). Both ADP tuning algorithms were parallelized and ran on the 16-core nodes of a cluster computer.

\begin{table}[]
\small
\caption{ADP performance on various types of systems, $\epsilon=0.01$}
\centering
\bgroup
\def\arraystretch{0.7}
\begin{tabular}{l|l|l|l|l|l|l|l}
\multicolumn{1}{c|}{$\pmb{\mu}$} & \multicolumn{1}{c|}{$\pmb{d}$} & \multicolumn{1}{c|}{$\pmb{TL}$} & \multicolumn{1}{c|}{\textbf{\begin{tabular}[c]{@{}c@{}}Performance\\ metric\end{tabular}}} & \multicolumn{1}{c|}{\textbf{CF}} & \multicolumn{1}{c|}{\textbf{Heuristic}} & \multicolumn{1}{c|}{\textbf{ADP GA}} & \multicolumn{1}{c}{\textbf{ADP TS}} \\ \hline
\multirow{12}{*}{0.1}            & \multirow{6}{*}{0.3}            & \multirow{3}{*}{20}              & Cost                                                                                       & 1.68                            & 1.26                                     & 1.35                                     & 1.40                                \\
                                 &                                 &                                  & FLAR                                                                                       & 97.0\%                           & 79.4\%                                  & 74.4\%                                 & 73.8\%                              \\
                                 &                                 &                                  & ART                                                                                        & 90.7                             & 65.1                                    & 78.5                                 & 84.6                                \\ \cline{3-8} 
                                 &                                 & \multirow{3}{*}{5}               & Cost                                                                                       & 0.49                             & 0.39                                    & 0.38                                 & 0.35                                \\
                                 &                                 &                                  & FLAR                                                                                       & 46.1\%                           & 36.7\%                                  & 34.3\%                               & 32.7\%                              \\
                                 &                                 &                                  & ART                                                                                        & 6.9                              & 5.7                                     & 5.8                                  & 6.1                                 \\ \cline{2-8} 
                                 & \multirow{6}{*}{2}              & \multirow{3}{*}{20}              & Cost                                                                                       & 1.07                             & 0.85                                    & 0.56                                 & 0.63                                \\
                                 &                                 &                                  & FLAR                                                                                       & 87.3\%                           & 60.2\%                                  & 40.7\%                               & 46.3\%                              \\
                                 &                                 &                                  & ART                                                                                        & 39.4                             & 42.4                                    & 32.2                                 & 33.2                                \\ \cline{3-8} 
                                 &                                 & \multirow{3}{*}{5}               & Cost                                                                                       & 0.15                             & 0.16                                    & 0.14                                 & 0.15                                \\
                                 &                                 &                                  & FLAR                                                                                       & 15.0\%                           & 15.7\%                                  & 13.1\%                               & 14.3\%                              \\
                                 &                                 &                                  & ART                                                                                        & 3.5                              & 3.3                                     & 4.5                                  & 4.6                                 \\ \hline
\multirow{12}{*}{0.5}            & \multirow{6}{*}{0.3}            & \multirow{3}{*}{20}              & Cost                                                                                       & 1.46                             & 1.58                                    & 0.85                                 & 0.96                                \\
                                 &                                 &                                  & FLAR                                                                                       & 95.6\%                           & 80.8\%                                  & 49.2\%                               & 58.6\%                              \\
                                 &                                 &                                  & ART                                                                                        & 69.7                             & 96.2                                    & 52.1                                 & 54.6                                \\ \cline{3-8} 
                                 &                                 & \multirow{3}{*}{5}               & Cost                                                                                       & 0.11                             & 0.13                                    & 0.09                                 & 0.10                                \\
                                 &                                 &                                  & FLAR                                                                                       & 10.4\%                           & 10.1\%                                  & 9.1\%                                & 9.3\%                               \\
                                 &                                 &                                  & ART                                                                                        & 3.5                              & 5.4                                     & 3.2                                  & 3.7                                 \\ \cline{2-8} 
                                 & \multirow{6}{*}{2}              & \multirow{3}{*}{20}              & Cost                                                                                       & 0.80                             & 0.88                                    & 0.30                                 & 0.31                                \\
                                 &                                 &                                  & FLAR                                                                                       & 68.2\%                           & 70.9\%                                  & 23.3\%                               & 24.2\%                              \\
                                 &                                 &                                  & ART                                                                                        & 29.2                             & 27.7                                    & 21.9                                 & 22.2                                \\ \cline{3-8} 
                                 &                                 & \multirow{3}{*}{5}               & Cost                                                                                       & 0.02                             & 0.02                                    & 0.02                                 & 0.02                                \\
                                 &                                 &                                  & FLAR                                                                                       & 2.0\%                            & 1.9\%                                   & 1.7\%                                & 1.8\%                               \\
                                 &                                 &                                  & ART                                                                                        & 2.7                              & 3.0                                     & 3.4                                  & 3.2                                
\end{tabular}
\egroup
\label{tab:num:agg01}
\end{table}

Table~\ref{tab:num:agg01} presents the obtained results. Note that our main objective was to minimize FLAR. Therefore, we chose the value of $\epsilon$ such that there is no notable decrease in ART compared to the other two policies. The ADP policy obtained with both tuning algorithms significantly outperforms both CF and Heuristic in terms of FLAR for all considered types of the systems. There is almost no increase in ART, and in some cases ART is also decreased by a large margin. The ADP policy performs especially good for the systems with larger distances and higher load, those where relocations contribute most to reducing response time to future incidents.

As already mentioned before, the ADP cost structure~\eqref{eq:costs} allows flexibility in terms of prioritizing FLAR against ART. If, for example, it is not important how large the waiting time of a broken machine is, given that it is over $TL$ time units, the $\epsilon$ parameter can be reduced to reflect that. Table~\ref{tab:num:agg001} shows the performance of the ADP policy for a subset of systems used in Table~\ref{tab:num:agg01}, with $\epsilon=0.001$.

\begin{table}[]
\small
\caption{ADP performance with $\epsilon=0.001$}
\centering
\bgroup
\def\arraystretch{0.7}
\begin{tabular}{l|l|l|l|l|l|l}
\multicolumn{1}{c|}{$\pmb{\mu}$} & \multicolumn{1}{c|}{$\pmb{d}$} & \multicolumn{1}{c|}{$\pmb{TL}$} & \multicolumn{1}{c|}{\textbf{\begin{tabular}[c]{@{}c@{}}Performance\\ metric\end{tabular}}} & \multicolumn{1}{c|}{\textbf{CF}} & \multicolumn{1}{c|}{\textbf{Heuristic}} & \multicolumn{1}{c}{\textbf{ADP GA}} \\ \hline
\multirow{12}{*}{0.1}            & \multirow{6}{*}{0.3}            & \multirow{3}{*}{20}              & Cost                                                                                       & 1.04                             & 0.84                                    & 0.49                                \\
                                 &                                 &                                  & FLAR                                                                                       & 97.0\%                           & 80.1\%                                  & 28.5\%                              \\
                                 &                                 &                                  & ART                                                                                        & 90.7                             & 61.1                                    & 225.1                               \\ \cline{3-7} 
                                 &                                 & \multirow{3}{*}{5}               & Cost                                                                                       & 0.46                             & 0.34                                    & 0.38                                \\
                                 &                                 &                                  & FLAR                                                                                       & 46.1\%                           & 34.3\%                                  & 33.8\%                              \\
                                 &                                 &                                  & ART                                                                                        & 6.9                              & 5.3                                     & 44.8                                \\ \cline{2-7} 
                                 & \multirow{6}{*}{2}              & \multirow{3}{*}{20}              & Cost                                                                                       & 0.89                             & 0.63                                    & 0.20                                \\
                                 &                                 &                                  & FLAR                                                                                       & 87.3\%                           & 59.9\%                                  & 9.3\%                               \\
                                 &                                 &                                  & ART                                                                                        & 39.4                             & 44.2                                    & 124.1                               \\ \cline{3-7} 
                                 &                                 & \multirow{3}{*}{5}               & Cost                                                                                       & 0.15                             & 0.15                                    & 0.12                                \\
                                 &                                 &                                  & FLAR                                                                                       & 15.0\%                           & 15.2\%                                  & 12.2\%                              \\
                                 &                                 &                                  & ART                                                                                        & 3.5                              & 3.3                                     & 4.4                                
\end{tabular}
\egroup
\label{tab:num:agg001}
\end{table}

For all the considered systems, both the genetic algorithm and the tabu search were able to find good solutions within a couple of days, and in some cases within a few hours. Note that the parameters $N^{GA}$ and $N^{TS}$ of the ADP tuning algorithms that affect both solution quality and computation time, as well as the number of performed iterations, were chosen such that the good policies are obtained for all systems within a reasonable amount of time. For any given system, the ADP policy can be further improved by increasing $N^{GA}$ and $N^{TS}$ and/or the number of iterations. The choice between the genetic algorithm and the tabu search depends on the computational resources available. In our experiments, when parallelized and ran on the 16-core nodes of a cluster computer, the genetic algorithm was able to find good quality solution about twice faster than the parallelized tabu search. However, when run sequentially, the tabu search was a few times faster.

\section{Conclusion}\label{sec:conclusion}
In this paper, we studied the problem of dynamic dispatching and relocation of service engineers. We considered the model introduced in~\cite{pechina2019} and developed an ADP approach to the problem. To that end, we proposed a number of basis functions and demonstrated in the computational experiments which of them are most important for finding a high performance policy. We also introduced two algorithms for tuning the coefficients of the basis functions in the ADP approach. We conducted extensive computational experiments, where we studied the parameters of the proposed ADP tuning algorithms, and compared the ADP approach against the two benchmark policies, closest-first policy and a heuristic algorithm that proved to perform well in~\cite{pechina2019}.

We showed that the ADP based policy outperforms both benchmarks for various types of systems, with most significant improvements for those with larger distances and higher loads. For the types of systems used in the study, we observed significant improvement over the best benchmark in terms of FLAR, without loss in ART. As it is computed offline and only once for each given type of system, the ADP policy is computationally tractable for real-life applications. Additionally, by modifying the cost structure with a single parameter $\epsilon$, it is possible to strike the desirable balance between the fraction of late arrivals and the average response time.
\newline

\noindent
\textbf{Acknowledgements}

This research was funded by an NWO grant, under contract number 438-15-506.

\bibliographystyle{plain}
\bibliography{references}

\begin{thebibliography}{10}

\bibitem{andersson2007}
Tobias Andersson and Peter V\"{a}rbrand.
\newblock Decision support tools for ambulance dispatch and relocation.
\newblock {\em Journal of the Operational Research Society}, 58(2):195--201,
  2007.

\bibitem{astaraky2015}
Davood Astaraky and Jonathan Patrick.
\newblock A simulation based approximate dynamic programming approach to
  multi-class, multi-resource surgical scheduling.
\newblock {\em European Journal of Operational Research}, 245(1):309--319,
  2015.

\bibitem{axsater2003}
Sven Axs{\"a}ter.
\newblock A new decision rule for lateral transshipments in inventory systems.
\newblock {\em Management Science}, 49(9):1168--1179, 2003.

\bibitem{bandara2014}
Damitha Bandara, Maria~E Mayorga, and Laura~A McLay.
\newblock Priority dispatching strategies for {EMS} systems.
\newblock {\em Journal of the Operational Research Society}, 65(4):572--587,
  2014.

\bibitem{belanger2019}
Val{\'e}rie B{\'e}langer, A~Ruiz, and Patrick Soriano.
\newblock Recent optimization models and trends in location, relocation, and
  dispatching of emergency medical vehicles.
\newblock {\em European Journal of Operational Research}, 272(1):1--23, 2019.

\bibitem{drent2018}
Collin Drent, Minou~Olde Keizer, and Geert-Jan van Houtum.
\newblock Dynamic dispatching and repositioning policies for service engineers.
\newblock Manuscript submitted for publication, 2018.

\bibitem{fang2013}
Jiarui Fang, Lei Zhao, Jan~C Fransoo, and Tom Van~Woensel.
\newblock Sourcing strategies in supply risk management: An approximate dynamic
  programming approach.
\newblock {\em Computers \& Operations Research}, 40(5):1371--1382, 2013.

\bibitem{gendreau1997}
Michel Gendreau, Gilbert Laporte, and Fr{\'e}d{\'e}ric Semet.
\newblock Solving an ambulance location model by tabu search.
\newblock {\em Location Science}, 5(2):75--88, 1997.

\bibitem{gendreau2001}
Michel Gendreau, Gilbert Laporte, and Fr{\'e}d{\'e}ric Semet.
\newblock A dynamic model and parallel tabu search heuristic for real-time
  ambulance relocation.
\newblock {\em Parallel Computing}, 27(12):1641--1653, 2001.

\bibitem{gendreau2006}
Michel Gendreau, Gilbert Laporte, and Fr{\'e}d{\'e}ric Semet.
\newblock The maximal expected coverage relocation problem for emergency
  vehicles.
\newblock {\em Journal of the Operational Research Society}, 57(1):22--28,
  2006.

\bibitem{glover1998}
Fred Glover and Manuel Laguna.
\newblock Tabu search.
\newblock In {\em Handbook of combinatorial optimization}, pages 2093--2229.
  Springer, 1998.

\bibitem{ingolfsson2013}
Armann Ingolfsson.
\newblock {EMS} planning and management.
\newblock In {\em Operations Research and Health Care Policy}, pages 105--128.
  Springer, 2013.

\bibitem{jagtenberg2015}
Caroline~J Jagtenberg, Sandjai Bhulai, and Robert~D van~der Mei.
\newblock An efficient heuristic for real-time ambulance redeployment.
\newblock {\em Operations Research for Health Care}, 4:27--35, 2015.

\bibitem{jagtenberg2017a}
Caroline~J Jagtenberg, Sandjai Bhulai, and Robert~D van~der Mei.
\newblock Dynamic ambulance dispatching: is the closest-idle policy always
  optimal?
\newblock {\em Health Care Management Science}, 20(4):517--531, 2017.

\bibitem{jagtenberg2017b}
Caroline~J Jagtenberg, Pieter~L van~den Berg, and Robert~D van~der Mei.
\newblock Benchmarking online dispatch algorithms for emergency medical
  services.
\newblock {\em European Journal of Operational Research}, 258(2):715--725,
  2017.

\bibitem{koudal2006}
Peter Koudal.
\newblock The service revolution in global manufacturing industries.
\newblock {\em Deloitte Research}, 2:1--22, 2006.

\bibitem{lam2007}
Shao-Wei Lam, Loo-Hay Lee, and Loon-Ching Tang.
\newblock An approximate dynamic programming approach for the empty container
  allocation problem.
\newblock {\em Transportation Research Part C: Emerging Technologies},
  15(4):265--277, 2007.

\bibitem{larson1974}
Richard~C. Larson.
\newblock A hypercube queuing model for facility location and redistricting in
  urban emergency services.
\newblock {\em Computers \& Operations Research}, 1(1):67--95, 1974.

\bibitem{larson1975}
Richard~C. Larson.
\newblock Approximating the performance of urban emergency service systems.
\newblock {\em Operations Research}, 23(5):845--868, 1975.

\bibitem{lee2007}
Young~Hae Lee, Jung~Woo Jung, and Young~Sang Jeon.
\newblock An effective lateral transshipment policy to improve service level in
  the supply chain.
\newblock {\em International Journal of Production Economics}, 106(1):115--126,
  2007.

\bibitem{maxwell2013}
Matthew~S Maxwell, Shane~G Henderson, and Huseyin Topaloglu.
\newblock Tuning approximate dynamic programming policies for ambulance
  redeployment via direct search.
\newblock {\em Stochastic Systems}, 3(2):322--361, 2013.

\bibitem{maxwell2010}
Matthew~S Maxwell, Mateo Restrepo, Shane~G Henderson, and Huseyin Topaloglu.
\newblock Approximate dynamic programming for ambulance redeployment.
\newblock {\em INFORMS Journal on Computing}, 22(2):266--281, 2010.

\bibitem{mclay2013b}
Laura~A McLay and Maria~E Mayorga.
\newblock A model for optimally dispatching ambulances to emergency calls with
  classification errors in patient priorities.
\newblock {\em IIE Transactions}, 45(1):1--24, 2013.

\bibitem{minner2003}
Stefan Minner, Edward~A Silver, and David~J Robb.
\newblock An improved heuristic for deciding on emergency transshipments.
\newblock {\em European Journal of Operational Research}, 148(2):384--400,
  2003.

\bibitem{murthy2004}
DNP Murthy, O~Solem, and T~Roren.
\newblock Product warranty logistics: Issues and challenges.
\newblock {\em European Journal of Operational Research}, 156(1):110--126,
  2004.

\bibitem{nasrollahzadeh2018}
Amir~Ali Nasrollahzadeh, Amin Khademi, and Maria~E Mayorga.
\newblock Real-time ambulance dispatching and relocation.
\newblock {\em Manufacturing \& Service Operations Management}, 20(3):467--480,
  2018.

\bibitem{pechina2019}
Anna Pechina, Dmitrii Usanov, Peter~M van~de Ven, and Rob~D van~der Mei.
\newblock Real-time dispatching and relocation of emergency service engineers.
\newblock {\em Manuscript submitted for publication}, 2019.

\bibitem{powell2007approximate}
Warren~B Powell.
\newblock {\em Approximate Dynamic Programming: Solving the Curses of
  Dimensionality}, volume 703.
\newblock John Wiley \& Sons, 2007.

\bibitem{schmid2012}
Verena Schmid.
\newblock Solving the dynamic ambulance relocation and dispatching problem
  using approximate dynamic programming.
\newblock {\em European Journal of Operational Research}, 219(3):611--621,
  2012.

\bibitem{schuetz2012}
Hans-Joerg Schuetz and Rainer Kolisch.
\newblock Approximate dynamic programming for capacity allocation in the
  service industry.
\newblock {\em European Journal of Operational Research}, 218(1):239--250,
  2012.

\bibitem{simao2009a}
Hugo Simao and Warren Powell.
\newblock Approximate dynamic programming for management of high-value spare
  parts.
\newblock {\em Journal of Manufacturing Technology Management}, 20(2):147--160,
  2009.

\bibitem{simao2009b}
Hugo~P Simao, Jeff Day, Abraham~P George, Ted Gifford, John Nienow, and
  Warren~B Powell.
\newblock An approximate dynamic programming algorithm for large-scale fleet
  management: A case application.
\newblock {\em Transportation Science}, 43(2):178--197, 2009.

\bibitem{tiemessen2013}
Harold~GH Tiemessen, Moritz Fleischmann, Geert-Jan van Houtum, Jo~AEE van
  Nunen, and Eleni Pratsini.
\newblock Dynamic demand fulfillment in spare parts networks with multiple
  customer classes.
\newblock {\em European Journal of Operational Research}, 228(2):367--380,
  2013.

\bibitem{barneveld2016b}
Thije van Barneveld.
\newblock The minimum expected penalty relocation problem for the computation
  of compliance tables for ambulance vehicles.
\newblock {\em INFORMS Journal on Computing}, 28(2):370--384, 2016.

\bibitem{barneveld2016a}
Thije van Barneveld, Sandjai Bhulai, and Robert~D van~der Mei.
\newblock The effect of ambulance relocations on the performance of ambulance
  service providers.
\newblock {\em European Journal of Operational Research}, 252(1):257--269,
  2016.

\bibitem{barneveld2017}
Thije van Barneveld, Sandjai Bhulai, and Robert~D van~der Mei.
\newblock A dynamic ambulance management model for rural areas.
\newblock {\em Health Care Management Science}, 20(2):165--186, 2017.

\bibitem{houtum2015}
Geert-Jan van Houtum and Bram Kranenburg.
\newblock {\em Spare Parts Inventory Control under System Availability
  Constraints}, volume 227.
\newblock Springer, 2015.

\bibitem{yu2010}
Xinjie Yu and Mitsuo Gen.
\newblock {\em Introduction to Evolutionary Algorithms}.
\newblock Springer Science \& Business Media, 2010.

\end{thebibliography}

\newpage
\appendix


\end{document}